\documentclass[11pt]{amsart}
\usepackage[T1]{fontenc}
\usepackage{lmodern,microtype}
\usepackage[margin=1in]{geometry}
\usepackage{amssymb,mathtools,booktabs,array,longtable}
\usepackage{xcolor,tikz,float,needspace}
\usepackage[colorlinks=true,linkcolor=blue!45!black,citecolor=blue!45!black,urlcolor=blue!45!black]{hyperref}
\hypersetup{pdftitle={The standard invariant of Petrescu's Hadamard matrix of order seven},pdfsubject={Hadamard matrices and relative commutants},pdflang={en}}
\newtheorem{theorem}{Theorem}[section]
\newtheorem{proposition}[theorem]{Proposition}
\newtheorem{lemma}[theorem]{Lemma}
\newtheorem{corollary}[theorem]{Corollary}
\theoremstyle{definition}

\theoremstyle{remark}
\newtheorem{remark}[theorem]{Remark}
\numberwithin{equation}{section}
\newcommand{\C}{\mathbb C}
\newcommand{\Q}{\mathbb Q}
\newcommand{\Z}{\mathbb Z}
\newcommand{\F}{\mathbb F}
\newcommand{\T}{\mathbb T}
\newcommand{\OO}{\mathcal O}
\newcommand{\RR}{\mathcal R}
\newcommand{\one}{\mathbf 1}
\newcommand{\Cat}{\operatorname{Cat}}
\newcommand{\TL}{\mathrm{TL}}
\DeclareMathOperator{\Tr}{Tr}
\DeclareMathOperator{\rank}{rank}

\title[Petrescu's Hadamard subfactors of order seven]{Invariants of 
Subfactors associated to \\ Petrescu's Hadamard matrices of order 7}
\author{Dietmar Bisch}
\address[Dietmar Bisch]{Department of Mathematics, Vanderbilt University, Nashville}
\thanks{Dietmar Bisch was supported in part by U.S. Army Research Office grant W911NF-23-1-0026.}
\author{Remus Nicoara}
\address[Remus Nicoara]{Department of Mathematics, University of Tennessee, Knoxville}
\date{}
\begin{document}

\begin{abstract}
Petrescu's sixth-root $7 \times 7$ Hadamard matrix defines an irreducible,
hyperfinite subfactor $N \subset M$ with Jones index 7. It is built
from the spin model commuting square obtained from the Hadamard matrix
by iterating Jones' basic construction. We prove that the standard
invariant of this subfactor is minimal, i.e. the higher relative
commutants consists of Temperley-Lieb algebras. Thus, the principal 
graphs of the subfactor are the Coxeter-Dynkin graph A$_\infty$. This 
solves a long-standing open problem. The main new idea is a congruence 
for the profile of the Hadamard matrix that restricts the order of rotation 
on the lowest-weight spaces in Jones planar algebra description of the
standard invariant. Together with a uniform
trace bound, this gives bounded annular multiplicities on a subsequence,
which forces the principal graph to have norm at most two.
We also show that the same result holds for most Hadamard matrices
in Petrescu's one-parameter family. Specifically, the associated
subfactors have Temperley-Lieb standard invariant at every
transcendental phase and at every algebraic phase with residue $\pm1$
at a prime above three. In each fixed degree the exceptional set is
finite and consists of algebraic numbers.

Our methods apply to several other families of Hadamard matrices
of order $q >4$ as well, and we show that the associated subfactors have
A$_\infty$ principal graphs. These include real Hadamard
matrices of order $q\equiv4\pmod8$, Hadamard matrices over
$\Q(i)$ of order $q\equiv2\pmod4$ and cube-root matrices of
order $q$ with $3\mid q$ and $9\nmid q$. For order six Hadamard
matrices, we
obtain the same result at every transcendental phase in
Di\c{t}\u{a}'s family through Haagerup's matrix and at every algebraically
independent pair of phases in Karlsson's two-parameter family.
\end{abstract}
\maketitle

\section{Introduction}
A Hadamard matrix of order $q$ is a matrix $H \in  M_q(\C)$ such
that its entries have modulus one and $HH^*=qI_q$. Then 
$U=q^{-1/2}H$ is of course a unitary matrix. We say that two
Hadamard matrices are equivalent if one can be obtained from the
other by permuting rows or columns, and by multiplying rows or
columns with scalars of absolute value one. Clearly, there is
a Hadamard matrix for every order, namely the Fourier matrix
(up to normalization). It is an intriguing question how many
equivalence classes of non-equivalent Hadamard matrices there
are of fixed order. If the order is a composite integer, it
is well-known, and not hard to see, that there are always
uncountably many non-equivalent Hadamard matrices of that order.
For $n=2$, $3$ and $5$, there is just the class of the Fourier
matrix \cite{H97}. If $n$ is a prime number $\ge 7$,
it was shown in \cite{HJ90,MW92}
that there is at least one other class of Hadamard matrices.
Petrescu showed in his thesis that there is in fact a 
one-parameter family of non-equivalent
Hadamard matrices of order $7$, $13$, $19$ etc. \cite{Pet97}.
It is not known what the answer is for order $11$.

The problem of studying complex Hadamard matrices was initiated
by Popa in \cite{P83} as they correspond to pairs of orthogonal MASAs in a
matrix algebra. This paper is at the origin of much of the
subsequent work on Hadamard matrices in the context of subfactor
theory.

Every $q \times q$ Hadamard matrix defines a so-called
{\it spin model commuting square} \cite[Chapter~6]{JS97} with which one
can build an irreducible, {\it hyperfinite} subfactor with Jones index $q$
\cite[Section~2]{N10}. With the above notation, the spin model commuting square associated to the Hadamard matrix $H$ is an inclusion of
finite dimensional matrix algebras of the form
\[
\begin{matrix}
UD_qU^*&\subset&M_q(\C)\\
\cup&&\cup\\
\C&\subset&D_q .
\end{matrix}
\]
where $D_q$ is the diagonal algebra. We fix the normalized matrix
trace on $M_q(\C)$ and iterate the basic construction 
horizontally (or vertically) to build a subfactor $N \subset M$ 
with index $q$ (see, for instance, \cite[Section~2]{N10},
\cite[Chapters~5 and~6]{JS97}, and \cite[Chapter~4]{GHJ89}).
Iterating the basic construction for $N \subset M$, we obtain
the Jones tower of II$_1$ factors
\[
N\subset M=M_0\subset M_1\subset M_2\subset\cdots
\]
We write $e_k\in M_k$ for the $k$th Jones projection, $k\geq1$, and set
\[
P_n=N'\cap M_{n-1}\quad(n\geq1),\qquad P_0=\C,
\]
$P_{n,-}=M'\cap M_n$ for $n\geq1$ and
$P_{0,-}=\C$. These centralizer algebras in the Jones tower 
associated to $N \subset M$ are called the {\it higher
relative commutants}. They capture {\it quantum symmetries} of
the subfactor as they are algebras of intertwiners of canonical
representations (bimodules) of the subfactor. They turn out to be 
intertwiner
spaces of certain unitary tensor categories of $N$-$N$ and
$M$-$M$ bimodules
and hence carry a group-like structure. The collection of higher
relative commutants, together with the data describing their
inclusions, traces, and Jones projections, is called the
{\it standard invariant} of the subfactor \cite[Section~1.3.7]{P94}.

Jones invented an algebraic-topological description of the standard
invariant and called it a {\it planar algebra} \cite[Theorem~4.2.1]{J21}. The higher relative
commutants are viewed as vector spaces $P_n$ and $P_{n,-}$ with 
black resp. white shading, and an action of the planar operad on the
collection of these spaces. In this picture, an important role is
played by the (two-click) rotation that is crucial for our
work as well.

If a subfactor is constructed from a commuting square, as in the
case of spin model subfactors associated to Hadamard matrices, it 
is known to be hard to compute the standard invariant, or, 
equivalently, the quantum symmetries the subfactor encodes. In fact,
there are only a few constructions where this was achieved, most
recently in \cite{BC25}. Since the algebra generated
by the Jones projections $e_k$, $1 \le k \le n$, is always a 
subalgebra of $N' \cap M_n$, it is particularly interesting to 
investigate if the subfactor has additional symmetries. If the
index is $>4$, these subalgebras are the Temperley-Lieb algebras
with parameter the square root of the Jones index. Thus, they 
represent the minimal quantum symmetries that every subfactor
has.

Popa proved in \cite{P93} that for every number $>4$ there is a
(non-hyperfinite) subfactor whose standard invariant is minimal.
This result was a breakthrough, as it established
that Temperley-Lieb algebras with any parameter $>2$ can be 
realized as quantum symmetries of a concrete subfactor. It is
unknown for which parameters this can be achieved by 
{\it hyperfinite} subfactors (see e.g. \cite{BC25}). 
In this paper, we obtain new examples.

It has been a long-standing open problem to compute the standard
invariant of Hadamard subfactors constructed from spin models.
Since Petrescu's one-parameter family of non-equivalent Hadamard
matrices of order $7$ is the first continuous family of prime
order, it is particularly interesting to study the associated 
subfactors. Experts have conjectured that their standard invariants
should be Temperley-Lieb, and partial progress was obtained
for instance
in \cite[Introduction and Section~4]{N06} and 
\cite[Example~3.3]{M21}. We show in this paper that the
Hadamard subfactors of index $7$ built from spin models based
on Petrescu's family have Temperley-Lieb standard invariant
for most parameters (see below for a precise statement).

\begin{theorem}\label{thm:main}
Let $N\subset M$ be the subfactor associated to Petrescu's matrix
$P_7$ of \cite[Eq.~(96)]{TZ06}. Then
\[
P_n=\TL_n(\sqrt7)\qquad(n\geq0).
\]
Thus, the principal graphs are $A_\infty$, and the dimensions of the 
higher relative commutants are the Catalan numbers:
\begin{align}
\dim(N'\cap M_k)&=\frac1{k+2}\binom{2k+2}{k+1}
&& (k\geq0),\label{eq:dimensions}\\
\dim(M'\cap M_k)&=\frac1{k+1}\binom{2k}{k}
&& (k\geq0).\label{eq:dual-dimensions}
\end{align}
\end{theorem}

We now describe the main ingredients of the proof. The first is the
observation that Petrescu's matrix $P_7$ has an equivalent form with
a particularly simple reduction modulo an Eisenstein prime. Let
$w=e^{\pi i/3}$ and put $\pi=1+w$. For the representative $H$ with
circulant core constructed in Section~\ref{sec:matrix}, reduction
modulo $\pi$ gives $-1$ on the diagonal and $1$ everywhere else.
Since $\Z[w]/(\pi)\cong\F_3$, we can write this as
\[
H\equiv I_7+J_7\pmod{\pi},
\]
where $I_7$ is the identity matrix and $J_7$ is the all-ones matrix.
This determines the profile coefficients of $H$ modulo $\pi$.

Using the cyclic coordinates of Section~\ref{sec:coordinates}, we
deduce that two-click rotation $\rho$ satisfies $\rho^3=I$ on $W_n$,
the common kernel of the cap maps in $P_n$. Since $\rho^n=I$ as well,
$\rho$ acts as the identity on $W_n$ whenever $3\nmid n$.

The other ingredient is the uniform estimate $|\Tr(\rho|P_n)|\leq7$ of
Section~\ref{sec:arithmetic}. Using Jones's annular decomposition~\cite{J01},
we show in Lemma~\ref{lem:odd-trace} that for every odd $n\geq3$,
$\Tr(\rho|P_n)=2+\Tr(\rho|W_n)$. Hence
$\dim W_n\leq5$ along an unbounded sequence of degrees. We use Jones's annular decomposition a second time through his formula
for $\dim P_n$ in terms of the annular multiplicities $\dim W_k$,
$k\leq n$ \cite[Theorem~6.3]{J01}. Since $\dim P_n$ also counts
the closed walks of length $2n$ based at the root of the principal
graph, we deduce in Proposition~\ref{prop:norm} that if the graph
had norm greater than two, then $\dim W_n$ would grow at least
exponentially with $n$. The bound above therefore forces the
principal graph to have norm at most two. Since $\dim P_2=2$ in our case, the graph is
$A_\infty$, and the Catalan dimensions give the theorem.

In Sections~\ref{sec:family} and~\ref{sec:algebraic-family}, we
investigate which other members of Petrescu's one-parameter family $H(z)$
give rise to subfactors with Temperley--Lieb standard invariant. We identify each relative commutant $P_n(z)$ with the kernel of a matrix
whose entries are Laurent polynomials in the parameter $z$ with
coefficients in $\Q(w)$. At each depth, the sixth-root
matrix supplies a nonvanishing minor of maximal rank. Thus, for each fixed $n$, we have $P_n(z)=\TL_n(\sqrt{7})$ for all
$z\in\T$ except possibly at finitely many algebraic values of $z$
(Theorem~\ref{thm:family}). The matrix congruence \eqref{eq:H-reduction} from Section~\ref{sec:matrix} extends
to any algebraic phase reducing to $1$ or $-1$ at a prime above
three (Theorem~\ref{thm:algebraic-family}). This includes all values
in $\Q(w)$ and all roots of unity of order dividing $2\cdot3^r$ for some $r\geq0$.

In Section~\ref{sec:extensions}, we apply the same methods to other
Hadamard subfactors of index $q>4$. For real matrices with $q\equiv4\pmod8$, Gaussian-rational
matrices with $q\equiv2\pmod4$, and cube-root matrices with $3\mid q$
and $9\nmid q$, the required congruences hold after division of the
profile by $4$, $2$, and $3$, respectively. Theorem~\ref{thm:ext-classes}
also treats algebraic representatives with local unit entries reducing
to $J_q-2I_q$ in odd characteristic $p$ when $q\equiv4\pmod p$.
Each congruence forces rotation to fix $W_n$ in odd degrees prime to
the residue characteristic. The resulting inclusions have standard
invariant $\TL(\sqrt q)$. In the real case Montgomery had obtained
$\dim P_2=2$ from the profile \cite[Proposition~3.2.9]{M22}.
The all-degree conclusion proves his conjecture for Paley type~II
matrices \cite[Conjecture~3.5]{M21}. Further examples include
Tao's matrix \cite{T04} and the matrices
$C+iI_q$ obtained from real symmetric conference matrices of the
indicated orders. By specialization we obtain $\TL(\sqrt6)$ for
every transcendental phase in Di\c{t}\u{a}'s affine family through
Haagerup's matrix and for algebraically independent phase pairs in
Karlsson's two-parameter family. Further algebraic parameters are given in
Lemma~\ref{lem:ext-haagerup} and Proposition~\ref{prop:ext-karlsson}.
Table~\ref{tab:ext-catalogue} lists the 844 representatives from the
Butson catalogue to which our criteria apply.

Throughout the paper we use the normalized traces on the Jones tower, write $\Tr$ for the
ordinary operator trace, and take inner products to be linear in the
second variable. 

\medskip
\noindent\textbf{AI Statement.}
OpenAI's ChatGPT assisted with the development and exposition of the
arguments and with checking references. The authors take responsibility
for the mathematical content of the paper.

\section{Petrescu's matrix and its profile}\label{sec:matrix}
Put
\[
w=e^{\pi i/3},\qquad w^2=w-1,\qquad
X=\{\infty\}\sqcup\Z/6\Z.
\]
With $\one=(1,\ldots,1)^{\mathsf T}\in\C^6$, consider the matrix
\begin{equation}\label{eq:H}
H=\begin{pmatrix}-1&\one^{\mathsf T}\\ \one&C\end{pmatrix},
\qquad C_{ij}=c_{j-i},\qquad
(c_0,\ldots,c_5)=(w,1,-w,w-1,-w,1).
\end{equation}
The core $C$ is a symmetric circulant. In terms of the permutation
matrix $R_{ij}=\one_{\{j=i+1\}}$ it is given by
\begin{equation}\label{eq:core}
C=wI+(R+R^{-1})-w(R^2+R^{-2})+(w-1)R^3.
\end{equation}

\begin{lemma}\label{lem:matrix}
The matrix \eqref{eq:H} is a sixth-root Hadamard matrix equivalent to
Petrescu's $P_7$ of \cite[Eq.~(96)]{TZ06}.
\end{lemma}
\begin{proof}
With indices starting at zero, take the rows and columns of $P_7$
in the respective orders
\[
(0,1,5,2,4,6,3),\qquad (0,2,5,1,3,6,4),
\]
then change the signs of the last six rows and the first column.
This gives \eqref{eq:H}, as the exponent matrices in
Appendix~\ref{app:equivalence} show. To check orthogonality directly,
the Fourier eigenvalues $\lambda_k=\sum_{s=0}^5c_sw^{ks}$ of $C$ are
\begin{equation}\label{eq:core-spectrum}
(\lambda_0,\ldots,\lambda_5)
=(1,\ 2+w,\ 3w-2,\ -1-2w,\ 3w-2,\ 2+w).
\end{equation}
Thus $C\one=\one$ and $CC^*=7I_6-J_6$, since $|\lambda_k|^2=7$
on $\one^\perp$. Block multiplication gives $HH^*=7I_7$.
\end{proof}

\subsection{The profile projection}\label{sec:profile}
We use the profile coefficients of \cite[Remark~2.5]{N10}, with
normalization
\begin{equation}\label{eq:profile}
S(a,b,c,d)=\sum_{y\in X}H_{ay}\overline{H_{by}}
                         \overline{H_{cy}}H_{dy},
\qquad p_{a,b}^{c,d}=\frac{S(a,b,c,d)}{49}.
\end{equation}
Regard the profile as a matrix on $\ell^2(X^2)$ by setting
\[
\Pi_{(a,c),(b,d)}=p_{a,b}^{c,d}.
\]

\begin{proposition}\label{prop:profile}
Set $\varepsilon_\infty=-1$, $\varepsilon_a=1$ for finite $a$, and
\[
g_\infty(r)=1,\qquad g_a(r)=c_{r-a}
\quad(a,r\in\Z/6\Z).
\]
The unnormalized profile is
\begin{equation}\label{eq:cyclic-profile}
S(a,b,c,d)=\varepsilon_a\varepsilon_b\varepsilon_c\varepsilon_d
 +\sum_{r=0}^5g_a(r)\overline{g_b(r)}\,
                    \overline{g_c(r)}g_d(r).
\end{equation}
The operator $\Pi$ is the rank-seven orthogonal projection
$\sum_{y\in X}v_yv_y^*$, where its range has orthonormal basis
\begin{equation}\label{eq:profile-vectors}
v_y(a,c)=\frac17H_{ay}\overline{H_{cy}},\qquad y\in X.
\end{equation}
\end{proposition}
\begin{proof}
The term indexed by $\infty$ in \eqref{eq:profile} gives the first
summand of \eqref{eq:cyclic-profile}. The other columns give the sum
over $r$. Column orthogonality gives
\[
\langle v_y,v_z\rangle
=\frac1{49}\left|\sum_a\overline{H_{ay}}H_{az}\right|^2
=\delta_{yz}.
\]
Formula \eqref{eq:profile} identifies the matrix entries of
$\sum_y v_yv_y^*$ with those of $\Pi$.
\end{proof}

Translation of the finite labels fixes $v_\infty$ and cyclically
permutes the remaining six vectors. Its Fourier decomposition gives
\begin{equation}\label{eq:fourier-profile}
u_k=\frac1{\sqrt6}\sum_{r=0}^5w^{-kr}v_r\quad(0\leq k\leq5),
\qquad
\Pi=v_\infty v_\infty^*+\sum_{k=0}^5u_ku_k^*.
\end{equation}
In the two-dimensional translation-fixed subspace of $\operatorname{ran}\Pi$ take the orthonormal vectors
\[
\mathfrak d=\frac{v_\infty+\sqrt6\,u_0}{\sqrt7},\qquad
\mathfrak t=\frac{\sqrt6\,v_\infty-u_0}{\sqrt7}.
\]
The identity $\sum_y v_y(a,c)=\delta_{ac}$ now yields
\begin{equation}\label{eq:split-profile}
\mathfrak d(a,c)=\frac{\delta_{ac}}{\sqrt7},\qquad
\mathfrak t(a,c)=\frac{\varepsilon_a\varepsilon_c-\delta_{ac}}{\sqrt{42}},
\qquad
\Pi=\mathfrak d\mathfrak d^*+\mathfrak t\mathfrak t^*
       +\sum_{k=1}^5u_ku_k^*.
\end{equation}
The diagonal pairs support the rank-one summand
$\mathfrak d\mathfrak d^*$. The complementary rank-six summand is
supported on ordered pairs of distinct labels.

\subsection{Reduction at the prime above three}
Let
\[
\OO=\Z[w],\qquad \pi=1+w.
\]
Since $w^2-w+1=0$, we have
\[
\OO/(\pi)
\cong\Z[t]/(t^2-t+1,t+1)\cong\F_3.
\]
Thus $\overline w$ and $w$ both reduce to $-1$. Also $\pi^2=3w$
and $\pi\overline\pi=3$. In \eqref{eq:H} the diagonal entries reduce
to $-1$ and the off-diagonal entries to $1$. Equivalently,
\begin{equation}\label{eq:H-reduction}
H\equiv I_7+J_7\pmod\pi.
\end{equation}

\begin{lemma}\label{lem:profile-congruence}
On pairs $a\ne c$, $b\ne d$, the profile satisfies
\begin{equation}\label{eq:S-reduction}
S(a,b,c,d)\equiv-1+\delta_{ab}+\delta_{ad}
                         +\delta_{cb}+\delta_{cd}\pmod\pi.
\end{equation}
\end{lemma}
\begin{proof}
Write $\mathbf e_a$ for the coordinate vectors of $\Z^X$ and $\one$
for the constant vector. The product of rows $a,c$ modulo $\pi$ is
$\one-2(\mathbf e_a+\mathbf e_c)$, since the two negative entries
occupy distinct positions. Conjugation fixes the residues, so
\begin{align*}
S(a,b,c,d)&\equiv
\bigl(\one-2(\mathbf e_a+\mathbf e_c)\bigr)\cdot
\bigl(\one-2(\mathbf e_b+\mathbf e_d)\bigr)\\
&=7-8+4(\delta_{ab}+\delta_{ad}+\delta_{cb}+\delta_{cd}).
\end{align*}
Reducing the last expression in $\F_3$ gives \eqref{eq:S-reduction}.
\end{proof}

\section{Cyclic coordinates for the relative commutants}\label{sec:coordinates}

\subsection{The angle operator and its fixed space}
Throughout this section $H$ is a complex Hadamard matrix of arbitrary
order $q$ with index set $X$. Write $D=D_q$ and denote the standard
matrix units of $M_q(\C)$ by $E_{ab}$.
The finite-dimensional tower of \cite[Propositions~2.2--2.3]{N10}
has algebras
\[
\mathcal P_{2r}=M_q(\C)^{\otimes(r+1)},\qquad
\mathcal P_{2r+1}=M_q(\C)^{\otimes(r+1)}\otimes D,
\qquad \mathcal P_{-1}=D.
\]
Let $U_i$ be the unitaries of \cite[Proposition~2.3]{N10}, with
$U_0=q^{-1/2}H$. Inside $\mathcal P_i$ put
\[
\mathcal A_i=D'\cap\mathcal P_i,\qquad
\mathcal Q_i=U_i\mathcal P_{i-1}U_i^*.
\]
Here $D$ acts on the first tensor factor. By \cite[Section~2]{N10},
\begin{equation}\label{eq:finite}
P_n\cong\mathcal A_{n-1}\cap\mathcal Q_{n-1}.
\end{equation}
For $n\geq2$ this is \cite[Proposition~2.4]{N10}. For $n=1$ the
intersection is $D\cap UDU^*=\C$.

We view both algebras in \eqref{eq:finite} as Hilbert subspaces of
$\mathcal P_{n-1}$ with the unnormalized Hilbert--Schmidt inner
product $\langle x,y\rangle=\Tr(x^*y)$. Scaling the trace in a
fixed degree leaves its orthogonal projections unchanged.
Since $D$ is maximal abelian in $M_q(\C)$,
\[
\mathcal A_{2r}=D\otimes M_q(\C)^{\otimes r},\qquad
\mathcal A_{2r+1}=D\otimes M_q(\C)^{\otimes r}\otimes D.
\]
In particular $\mathcal A_{n-1}$ and $\mathcal Q_{n-1}$ both have
dimension $q^n$. Let $\mathcal H_n=\ell^2(X^n)$ with counting measure.
An orthonormal basis in each algebra gives coordinate isometries
\[
\alpha_n:\mathcal H_n\longrightarrow\mathcal A_{n-1},\qquad
\beta_n:\mathcal H_n\longrightarrow\mathcal Q_{n-1}.
\]
Write $E_{\mathcal A},E_{\mathcal Q}$ for the trace-preserving
conditional expectations onto these algebras. As operators on the
ambient Hilbert space they are the range projections of $\alpha_n$
and $\beta_n$. Thus $C_n=\alpha_n^*\beta_n$ satisfies
\begin{equation}\label{eq:two-projections}
T_n:=C_nC_n^*=\alpha_n^*E_{\mathcal Q}\alpha_n,
\qquad
\alpha_nT_n=E_{\mathcal A}E_{\mathcal Q}\alpha_n.
\end{equation}
For $x=\alpha_nf$ one has
\begin{equation}\label{eq:distance}
\langle f,(I-T_n)f\rangle
=\|x\|_{\rm HS}^2-\|E_{\mathcal Q}x\|_{\rm HS}^2
=\|x-E_{\mathcal Q}x\|_{\rm HS}^2.
\end{equation}
It follows that $T_n$ is a positive contraction and that $\alpha_n$
gives an isometry from $\ker(I-T_n)$ onto the intersection in
\eqref{eq:finite} for the Hilbert--Schmidt inner product.
This is the angle-operator description of flatness
\cite[Theorem~2.11.8 and the ensuing spin-model remark]{J21}.
See also \cite[Section~1]{M21} for this formulation and \cite{M23}
for the commuting-square planar algebra. In cyclic matrix-unit bases we obtain
\begin{align}
C_n(\mathbf a,\mathbf z)
 &=q^{-(n+1)/2}\prod_{j=1}^n
 H_{a_jz_j}\overline{H_{a_{j+1}z_j}},\label{eq:Cn}\\
T_n(\mathbf a,\mathbf b)
 &=q^{-n-1}\prod_{j=1}^n
 S(a_j,b_j,a_{j+1},b_{j+1}).\label{eq:Tn}
\end{align}
All word indices are cyclic. Summing the first formula separately
in each coordinate of $\mathbf z$ gives the second,
\[
\begin{aligned}
(C_nC_n^*)(\mathbf a,\mathbf b)
&=q^{-n-1}\prod_{j=1}^n
\left(\sum_{z_j\in X}H_{a_jz_j}\overline{H_{b_jz_j}}
\overline{H_{a_{j+1}z_j}}H_{b_{j+1}z_j}\right)\\
&=T_n(\mathbf a,\mathbf b).
\end{aligned}
\]

\subsection{Matrix units and the flatness condition}
We recall the planar algebra realization of \eqref{eq:finite}.
Jones's spin planar algebra has modulus $\sqrt q$ and its operations
are given by summation over region labels
\cite[Example~2.8 and Proposition~2.8.1]{J21}. For a biunitary,
the flat elements form a planar $*$-subalgebra
\cite[Theorem~2.11.3]{J21}. The crossed-string and shift formulations
of flatness are given in \cite[Propositions~2.11.4 and~2.11.6]{J21}.
For a Hadamard matrix this subalgebra is a spherical $C^*$-planar
algebra, and its flatness equations are those of the commuting-square
relative commutants \cite[Proposition~2.11.12 and the paragraph
following it]{J21}.

Choose the shading whose distinguished region is a leaf of the star.
Each leaf $a\in X$ has weight one and the central vertex $\bullet$
has weight $\sqrt q$. A positive $n$-box basis element is a loop
\[
a_1,\bullet,a_2,\bullet,\ldots,a_n,\bullet,a_1.
\]
Its two halves, read from the distinguished leaf to the opposite
vertex, are the two paths of a matrix unit. According to the parity
of $n$ they give the algebras
\[
D\otimes M_q^{\otimes r}\quad(n=2r+1),\qquad
D\otimes M_q^{\otimes r}\otimes D\quad(n=2r+2).
\]
These identify with $\mathcal A_{n-1}$ in \eqref{eq:finite}.

\begin{lemma}\label{lem:coordinates}
In the path basis the flatness condition takes the form
\[
x=U_{n-1}\widetilde y\,U_{n-1}^*,\qquad
x\in\mathcal A_{n-1},\quad y\in\mathcal P_{n-2},
\]
where $\widetilde y$ denotes the image of $y$ under the tower
inclusion. The cyclic matrix units define isometries $\alpha_n,\beta_n$
onto $\mathcal A_{n-1},\mathcal Q_{n-1}$ with
$\alpha_n^*\beta_n=C_n$ as in \eqref{eq:Cn}.
\end{lemma}

\begin{proof}
Give a crossing weight $H_{xy}$ and its reverse weight
$\overline{H_{xy}}$, with change of basis $U=q^{-1/2}H$. Put
\[
D_U=\sum_{x,y}\overline{H_{xy}}E_{yy}\otimes E_{xx},\qquad
G=(U\otimes I)D_U.
\]
Then
\[
G_{(a,t),(b,s)}=\delta_{ts}\,q^{-1/2}H_{ab}\overline{H_{tb}}.
\]
On consecutive tensor factors the recursion of
\cite[Proposition~2.3]{N10} is
\begin{align*}
U_{2r+1}&=G_{0,1}G_{1,2}\cdots G_{r,r+1},\\
U_{2r}&=U_{2r-1}(I^{\otimes r}\otimes U)\quad(r\geq1).
\end{align*}
All intermediate indices are fixed by the deltas in $G$. Multiplication
in the indicated order gives
\begin{align}
(U_{2r})_{\mathbf x,\mathbf y}
 &=q^{-(r+1)/2}\prod_{j=0}^r H_{x_jy_j}
       \prod_{j=0}^{r-1}\overline{H_{x_{j+1}y_j}},\label{eq:strip-even}\\
(U_{2r+1})_{\mathbf x,\mathbf y}
 &=\delta_{x_{r+1},y_{r+1}}q^{-(r+1)/2}
       \prod_{j=0}^r H_{x_jy_j}\overline{H_{x_{j+1}y_j}}.
       \label{eq:strip-odd}
\end{align}
The strip on the row path contributes $U_i$. Reflection gives its
adjoint on the column path. Thus the coefficient on the exterior
paths $\mathbf x,\mathbf x'$ is
\begin{equation}\label{eq:flat-cut-entry}
\sum_{\mathbf y,\mathbf y'}
(U_i)_{\mathbf x,\mathbf y}\,
\widetilde y_{\mathbf y,\mathbf y'}\,
\overline{(U_i)_{\mathbf x',\mathbf y'}},\qquad i=n-1.
\end{equation}
This is the $(\mathbf x,\mathbf x')$ entry of
$U_i\widetilde yU_i^*$, to be equated with the entry of
$x\in\mathcal A_i$ on the other side of the flatness relation
(Figure~\ref{fig:flat-cut}).

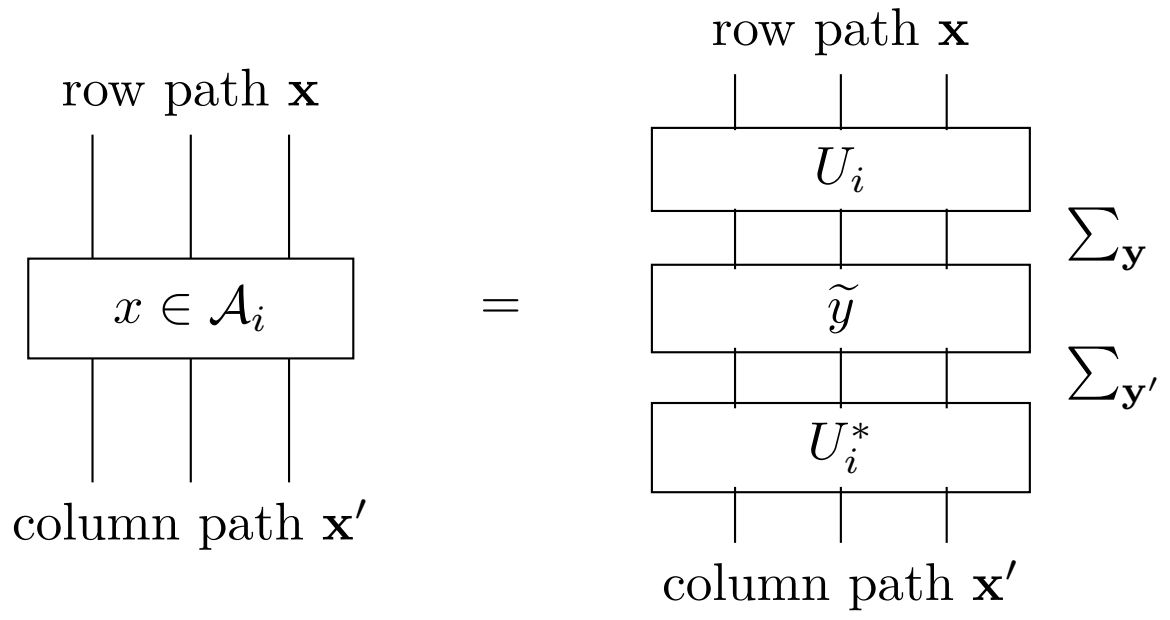
\begin{figure}[H]
\centering
\begin{tikzpicture}[x=1cm,y=1cm,>=stealth, every node/.style={font=\small}]
\node[draw,minimum width=2.15cm,minimum height=.66cm] (x) at (0,0) {$x\in\mathcal A_i$};
\foreach \t in {-.65,0,.65}{\draw (\t,.33)--(\t,1.15);\draw (\t,-.33)--(\t,-1.15);}
\node at (0,1.43) {row path $\mathbf x$};
\node at (0,-1.43) {column path $\mathbf x'$};
\node at (2.05,0) {$=$};
\node[draw,minimum width=2.5cm,minimum height=.52cm] (u) at (4.3,.92) {$U_i$};
\node[draw,minimum width=2.5cm,minimum height=.52cm] (y) at (4.3,0) {$\widetilde y$};
\node[draw,minimum width=2.5cm,minimum height=.52cm] (v) at (4.3,-.92) {$U_i^*$};
\foreach \t in {3.6,4.3,5}{\draw (\t,1.18)--(\t,1.55);\draw (\t,.66)--(\t,.26);\draw (\t,-.26)--(\t,-.66);\draw (\t,-1.18)--(\t,-1.55);}
\node[right] at (5.65,.46) {$\sum_{\mathbf y}$};
\node[right] at (5.65,-.46) {$\sum_{\mathbf y'}$};
\node at (4.3,1.83) {row path $\mathbf x$};
\node at (4.3,-1.83) {column path $\mathbf x'$};
\end{tikzpicture}
\caption{The matrix form of the flatness equation. Each group of lines
represents a tensor path. The boxes are the matrices evaluated in
\eqref{eq:strip-even} and \eqref{eq:strip-odd}, and the sums are the
matrix products in \eqref{eq:flat-cut-entry}.}
\label{fig:flat-cut}
\end{figure}

Suppose first that $i=n-1=2r$. The first path basis
is given by
\[
\alpha_{2r+1}e_{(a_0,\ldots,a_r,b_r,\ldots,b_1)}
=E_{a_0a_0}\otimes E_{a_1b_1}\otimes\cdots\otimes E_{a_rb_r}.
\]
For the second, a word $\mathbf z=(c_0,\ldots,c_{r-1},s,
d_{r-1},\ldots,d_0)$ gives the matrix unit
\[
Y=E_{c_0d_0}\otimes\cdots\otimes E_{c_{r-1}d_{r-1}}\otimes E_{ss},
\qquad \beta_{2r+1}e_{\mathbf z}=U_{2r}YU_{2r}^*.
\]
The inner product of these two elements is
\[
(U_{2r})_{(a_0,\ldots,a_r),(c_0,\ldots,c_{r-1},s)}
\overline{(U_{2r})_{(a_0,b_1,\ldots,b_r),(d_0,\ldots,d_{r-1},s)}}.
\]
By \eqref{eq:strip-even} the first factor reads off the row path
and the conjugated factor reads off the return path in reverse order.
Joining them at $s$ and $a_0$ gives the cyclic product \eqref{eq:Cn}.
At $r=0$ the coefficient is
$\langle E_{aa},UE_{zz}U^*\rangle=1/q$.

For $i=n-1=2r+1$ and $\mathbf a=(a_0,\ldots,a_r,t,b_r,\ldots,b_1)$
we have
\[
\alpha_{2r+2}e_{\mathbf a}
=E_{a_0a_0}\otimes E_{a_1b_1}\otimes\cdots
 \otimes E_{a_rb_r}\otimes E_{tt}.
\]
In the second basis the word $(c_0,\ldots,c_r,d_r,\ldots,d_0)$
corresponds to $q^{-1/2}E_{\mathbf c,\mathbf d}\otimes I$ conjugated
by $U_{2r+1}$. Indeed the embedded matrix unit before normalization
has squared norm $q$. The inner product is now
\[
q^{-1/2}(U_{2r+1})_{(a_0,\ldots,a_r,t),(c_0,\ldots,c_r,t)}
\overline{(U_{2r+1})_{(a_0,b_1,\ldots,b_r,t),(d_0,\ldots,d_r,t)}}.
\]
Formula~\eqref{eq:strip-odd} gives the same cyclic product, with
normalizing factor $q^{-1/2}q^{-(r+1)}=q^{-(n+1)/2}$.
\end{proof}

\subsection{Rotation and the common cap kernel}
Let $\rho$ be rotation by two consecutive boundary intervals. For
$n\geq1$ denote by $W_n$ the intersection in $P_n$ of the kernels
of all elementary caps joining adjacent boundary points.

\begin{proposition}\label{prop:transfer}
For $n\geq1$ one has $0\leq T_n=C_nC_n^*\leq I$ and
\begin{equation}\label{eq:fixed}
P_n\cong\mathcal F_n:=\ker(I-T_n).
\end{equation}
Under the identification $P_n\cong F_n$ in~\eqref{eq:fixed}, two-click
rotation corresponds, for one choice of orientation, to the restriction
to $F_n$ of the cyclic shift on $\ell^2(X^n)$ given by
\[
e_{(a_1,\ldots,a_n)}\longmapsto e_{(a_2,\ldots,a_n,a_1)}.
\]For $n\geq2$ one has
\begin{equation}\label{eq:W}
W_n=\mathcal F_n\cap\mathcal U_n,
\end{equation}
where $f\in\mathcal U_n$ precisely when
\[
f(a_1,\ldots,a_n)=0\quad\text{if }a_j=a_{j+1}\text{ for some }j,
\]
and
\[
\sum_{x\in X}f(a_1,\ldots,a_{j-1},x,a_{j+1},\ldots,a_n)=0
\]
for every $j$ and every choice of the remaining labels. The indices
are cyclic. In degree one, $W_1=0$.
\end{proposition}

\begin{proof}
The operator assertions follow from \eqref{eq:two-projections},
\eqref{eq:distance} and Lemma~\ref{lem:coordinates}. To compute
rotation, move the distinguished interval to the next leaf region.
The initial and final weights are both one. The resulting word begins
at $a_2$ and ends at $a_1$, as in Figure~\ref{fig:word-rotation}.
On coefficient functions the action is consequently
\[
(\rho f)(a_1,\ldots,a_n)=f(a_n,a_1,\ldots,a_{n-1}).
\]

\begin{figure}[H]
\centering
\begin{tikzpicture}[scale=.83,>=stealth,every node/.style={font=\small}]
\foreach \cx/\shift in {0/0,5/1}{
 \begin{scope}[shift={(\cx,0)}]
 \draw (0,0) circle (1.38);
 \foreach \angle in {0,45,90,135,180,225,270,315}
  {\fill (\angle:1.38) circle (.045);}
 \foreach \angle/\lab in {67.5/a_1,157.5/a_2,247.5/a_3,337.5/a_4}
  {\node at (\angle:1.02) {$\lab$};}
 \foreach \angle in {22.5,112.5,202.5,292.5}
  {\node at (\angle:1.03) {$\bullet$};}
 \ifnum\shift=0
  \node at (67.5:1.7) {$*$};
  \node at (0,-1.95) {$(a_1,a_2,a_3,a_4)$};
 \else
  \node at (157.5:1.7) {$*$};
  \node at (0,-1.95) {$(a_2,a_3,a_4,a_1)$};
 \fi
 \end{scope}
}
\draw[->,thick] (1.85,0)--(3.15,0) node[midway,above] {$\rho$};
\node[align=center,font=\footnotesize] at (2.5,-.62) {two\\endpoints};
\end{tikzpicture}
\caption{Rotation at level four. The points on each circle are the
string endpoints, and the labels belong to the intervening regions.
Moving the distinguished interval by two boundary points changes
the order of the leaf labels from $(a_1,a_2,a_3,a_4)$ to
$(a_2,a_3,a_4,a_1)$.}
\label{fig:word-rotation}
\end{figure}
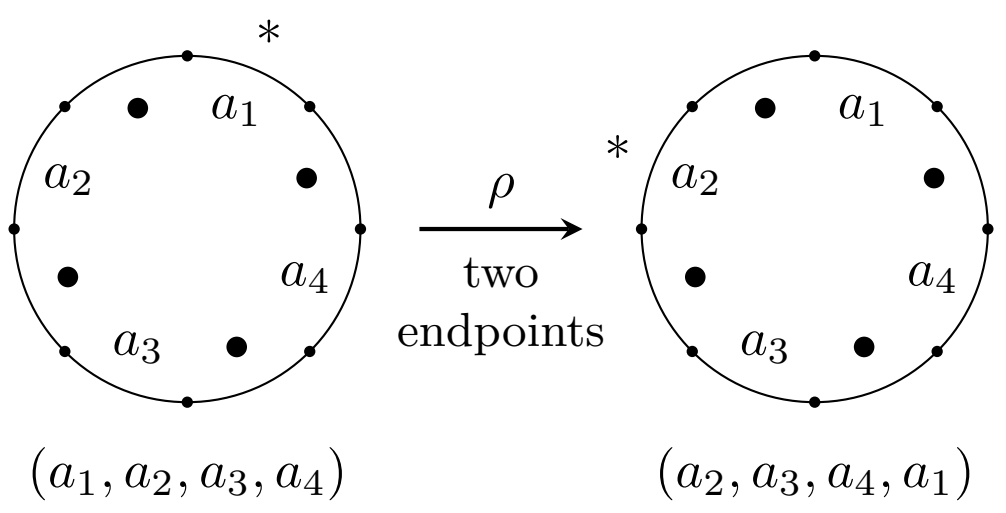

A cap either removes a leaf region, summing its label, or removes
a central region, identifying the two adjacent leaves
(Figure~\ref{fig:two-caps}). Caps meeting the distinguished interval
are obtained from these by rotation.

\begin{figure}[H]
\centering
\begin{tikzpicture}[x=1cm,y=1cm,every node/.style={font=\small}]
\begin{scope}
 \draw[thick] (0,0)--(3.1,0);
 \node[below] at (1.55,-.03) {input boundary};
 \fill[gray!16] (.85,0)..controls(.85,1.45)and(2.25,1.45)..(2.25,0)--cycle;
 \draw[thick] (.85,0)..controls(.85,1.45)and(2.25,1.45)..(2.25,0);
 \fill (.85,0) circle(.035);\fill(2.25,0) circle(.035);
 \node at (1.55,.43) {$x\in X$};
 \node at (.36,.37) {$\bullet$};\node at (2.74,.37) {$\bullet$};
 \node[align=center] at (1.55,1.65) {closed leaf region\\sum over $x$};
 \node at (1.55,-.92) {$q^{-1/2}\sum_x f(\ldots,x,\ldots)$};
\end{scope}
\begin{scope}[shift={(5.1,0)}]
 \draw[thick] (0,0)--(3.1,0);
 \node[below] at (1.55,-.03) {input boundary};
 \draw[thick] (.85,0)..controls(.85,1.45)and(2.25,1.45)..(2.25,0);
 \fill (.85,0) circle(.035);\fill(2.25,0) circle(.035);
 \node at (1.55,.43) {$\bullet$};
 \node at (.36,.37) {$a$};\node at (2.74,.37) {$b$};
 \draw[densely dotted] (.36,.62)..controls(.36,1.3)and(2.74,1.3)..(2.74,.62);
 \node[align=center] at (1.55,1.65) {one exterior leaf region\\requires $a=b$};
 \node at (1.55,-.92) {$f(\ldots,a,a,\ldots)$};
\end{scope}
\end{tikzpicture}
\caption{The two local cap contractions. The dotted arc on the
right indicates a path within the exterior region, not an additional
string. The displayed formulas use the counting normalization;
the geometric caps are nonzero scalar multiples of them. Other
boundary labels pass through unchanged.}
\label{fig:two-caps}
\end{figure}
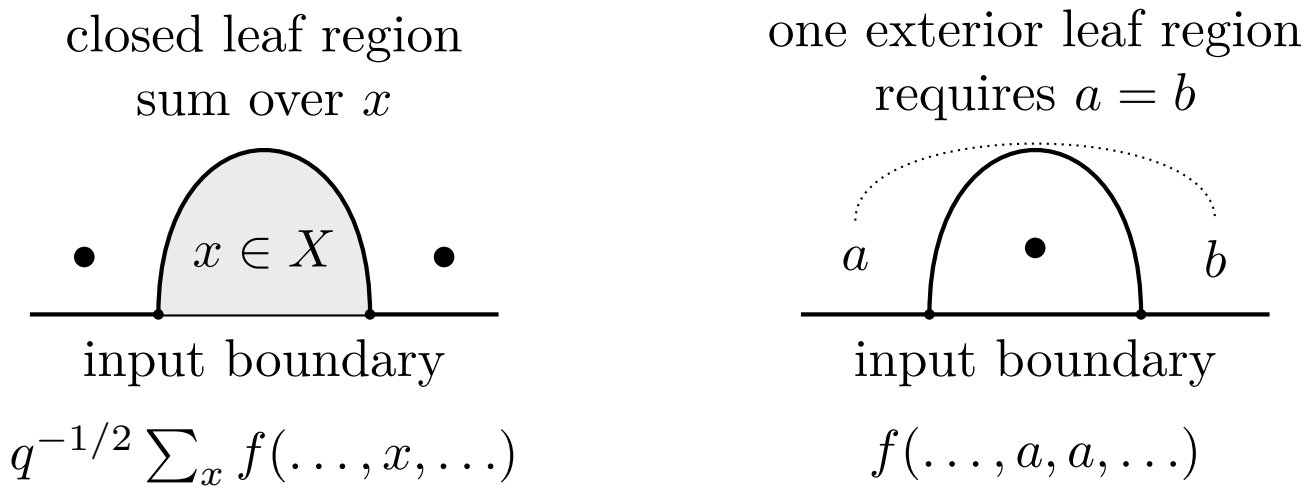

For $n\geq 2$, the corresponding normalized cups are isometries
$I_j,D_j:\mathcal H_{n-1}\to\mathcal H_n$. For $1\leq j\leq n$ put
\[
(I_jh)(a_1,\ldots,a_n)=q^{-1/2}
h(a_1,\ldots,\widehat{a_j},\ldots,a_n).
\]
For $j<n$ put
\[
(D_jh)(a_1,\ldots,a_n)=\one_{\{a_j=a_{j+1}\}}
h(a_1,\ldots,a_j,a_{j+2},\ldots,a_n).
\]
At the distinguished interval the formula is
\[
(D_nh)(a_1,\ldots,a_n)=\one_{\{a_n=a_1\}}
h(a_1,\ldots,a_{n-1}).
\]
Taking adjoints for counting measure gives
\begin{align*}
(I_j^*f)(b_1,\ldots,b_{n-1})
 &=q^{-1/2}\sum_{x\in X}
 f(b_1,\ldots,b_{j-1},x,b_j,\ldots,b_{n-1}),\\
(D_j^*f)(b_1,\ldots,b_{n-1})
 &=f(b_1,\ldots,b_j,b_j,b_{j+1},\ldots,b_{n-1})\quad(j<n).
\end{align*}
In particular the cup-cap projections $A_j=I_jI_j^*$ and $E_j=D_jD_j^*$
are given by
\[
(A_jf)(\mathbf a)=q^{-1}\sum_x
 f(a_1,\ldots,a_{j-1},x,a_{j+1},\ldots,a_n),\qquad
(E_jf)(\mathbf a)=\one_{\{a_j=a_{j+1}\}}f(\mathbf a).
\]
All leaves have the same weight, so a geometric cap differs from
the corresponding adjoint only by a nonzero scalar. The normalized
trace and counting inner products also differ by a common factor
at each degree. Neither normalization affects the kernels.

One may see the restriction to the flat spaces directly. Orthogonality
of the rows and columns of $H$ gives
\[
C_nA_j=E_jC_n,\qquad C_nE_j=A_{j+1}C_n,
\]
with cyclic subscripts. It follows that $T_n$ commutes with the
$A_j$ and $E_j$. More precisely, the identities
\[
S(a,b,c,b)=q\delta_{ac},\qquad
\sum_x S(a,b,c,x)S(c,x,d,e)=qS(a,b,d,e)
\]
give
\[
T_nD_j=D_jT_{n-1},\qquad T_nI_j=I_jT_{n-1}.
\]
The first uses the repeated label at a $D_j$ cup. In the other,
summing the two profile factors over their common label contracts
them by $\sum_xH_{xy}\overline{H_{xz}}=q\delta_{yz}$. Taking adjoints
of the two intertwining identities shows that both caps map
$\mathcal F_n$ into $\mathcal F_{n-1}$. Their kernels in the ambient
space are
\[
\ker I_j^* = \{f:\text{the $j$th coordinate sum is zero}\},\qquad
\ker D_j^* = \{f:f|_{\{a_j=a_{j+1}\}}=0\}.
\]
Intersecting with $\mathcal F_n$ proves \eqref{eq:W}. Finally,
$T_1=J_q/q$ has the constants as its fixed space. The cap is nonzero
on the identity, whence $W_1=0$.
\end{proof}

\begin{lemma}\label{lem:low}
For the Petrescu inclusion, $\dim P_1=1$ and $\dim P_2=2$.
\end{lemma}
\begin{proof}
The first assertion follows from $T_1=J_7/7$. We give the profile
component calculation for $P_2$ in the coordinates above
\cite[Theorem~2.11.15]{J21}. In $T_2$ the diagonal
pairs form a block $J_7/7$. On its complement the entries are
\[
T_2((a,c),(b,d))=\frac{|S(a,b,c,d)|^2}{7^3}.
\]
The restriction $\mathsf S=49\Pi$ to this complement satisfies
$\mathsf S^2=49\mathsf S$ by Proposition~\ref{prop:profile}.
Its diagonal entries are seven, so
\[
\sum_{b\ne d}|S(a,b,c,d)|^2=49\cdot7=7^3.
\]
Thus the second block is stochastic. By \eqref{eq:S-reduction} its
support joins disjoint pairs as well as the two orientations of any
pair. Any two pairs have a third pair disjoint from both, so the
block is irreducible. Its fixed space is one-dimensional by
Perron--Frobenius. Together the two blocks give $\dim P_2=2$.
\end{proof}

\section{Rotation on the relative commutants}\label{sec:arithmetic}
\subsection{Rotation on the lowest-weight space}
Write
\[
\Omega_n=\{(a_1,\ldots,a_n)\in X^n:a_j\ne a_{j+1}\text{ for all }j\},
\]
where $a_{n+1}=a_1$. By \eqref{eq:W} every element of $W_n$ is
supported on $\Omega_n$. On this set form the matrix
\[
K_n(\mathbf a,\mathbf b)=\prod_{j=1}^n
S(a_j,b_j,a_{j+1},b_{j+1}).
\]
It has coefficients in $\OO=\Z[w]$. Clearing the denominator in
\eqref{eq:Tn} gives the flatness equation
\begin{equation}\label{eq:integral-flat}
K_nf=7^{n+1}f\qquad(f\in W_n).
\end{equation}
We work over the local ring
\[
\RR=\OO_{(\pi)}
=\{a/b:a,b\in\OO,\ b\notin(\pi)\},\qquad \pi=1+w.
\]
The norm $|a+bw|^2=a^2+ab+b^2$ is Euclidean on $\OO$, so $\RR$
is a discrete valuation ring with uniformizer $\pi$ and residue
field $\F_3$.

\begin{lemma}\label{lem:cyclic-expansion}
For $n\geq2$ and $f\in\mathcal U_n$ with entries in $\RR$
(in particular, in $\OO$), one has on $\Omega_n$
\begin{equation}\label{eq:cyclic-congruence}
K_nf\equiv(2I+\rho+\rho^{-1})f\pmod\pi.
\end{equation}
\end{lemma}
\begin{proof}
Fix $\mathbf a\in\Omega_n$. In the sum defining $(K_nf)(\mathbf a)$,
replace each factor by the right-hand side of \eqref{eq:S-reduction}.
After this replacement we may sum over all $\mathbf b$, extending
$f$ by zero off $\Omega_n$. A term independent of some $b_j$ has
zero sum by the cap equation in that coordinate. Each nonconstant
factor involves just one input label. Thus a surviving term uses
one delta in each factor and each input label exactly once.

Let $\epsilon_j$ be zero if the $j$th factor uses $b_j$ and one if
it uses $b_{j+1}$. The number of occurrences of $b_j$ is
$1-\epsilon_j+\epsilon_{j-1}$. Hence all $\epsilon_j$ coincide.
For either choice the selected inputs are now equated with $a_j$
or $a_{j+1}$. A change from the second choice in one factor to
the first in the next gives equal adjacent inputs and hence a zero
value of $f$. Every nonconstant cyclic choice has such a change.
There remain two identity terms and one term for each cyclic shift.
When $n=2$ the two shifts agree and both are counted.
\end{proof}

\begin{lemma}\label{lem:local-torsion}
Let $B$ be a finite-order automorphism of a finite free $\RR$-module
$L$. If $(B-I)^2L\subseteq\pi L$, then $B^3=I$.
\end{lemma}
\begin{proof}
Extend the valuation of $\Q(w)$ to its algebraic closure, with
$v(\pi)=1$. Then $v(3)=2$. Since $(B-I)^2/\pi$ preserves $L$,
each eigenvalue $\zeta$ of $B$ satisfies
\begin{equation}\label{eq:valuation}
v(\zeta-1)\geq\tfrac12.
\end{equation}
Write the order of $\zeta$ as $3^km$ with $3\nmid m$. The root
$\eta=\zeta^{3^k}$ has order $m$ and residue one. If $m>1$, the
relation $\sum_{j=0}^{m-1}\eta^j=0$ would give $m=0$ in the residue
field. Thus $m=1$.

For $k\geq1$ the factors $1-\zeta^a$ with $1\leq a<3^k$ and
$3\nmid a$ all have valuation $v(1-\zeta)$, since
$(1-\zeta^a)/(1-\zeta)$ reduces to $a$. Their product is
$\Phi_{3^k}(1)=3$. It follows that
\[
v(1-\zeta)=\frac{2}{\varphi(3^k)}=\frac1{3^{k-1}}.
\]
Comparison with \eqref{eq:valuation} gives $k\leq1$.
Since $B$ is diagonalizable in characteristic zero, $B^3=I$.
\end{proof}

\begin{theorem}\label{thm:rotation}
The restriction of two-click rotation to $W_n$, $n\geq2$, satisfies
\begin{equation}\label{eq:rotation-order}
\rho^3=I\quad\text{on }W_n.
\end{equation}
In particular, $\rho|_{W_n}=I$ if $3\nmid n$.
\end{theorem}
\begin{proof}
Let $W_{n,K}$ be the solution over $K=\Q(w)$ of the equations
\eqref{eq:fixed} and \eqref{eq:W}. They give
$W_n=W_{n,K}\otimes_K\C$. Consider the lattice
\[
L_n=W_{n,K}\cap\RR^{\Omega_n}.
\]
It is free over $\RR$ and spans $W_{n,K}$, as one sees by clearing
denominators in a basis. It is stable under $\rho$ and saturated
in $\RR^{\Omega_n}$, since division of a solution by $\pi$ still
satisfies the defining linear equations.

For $f\in L_n$, equations \eqref{eq:integral-flat} and
\eqref{eq:cyclic-congruence}, with $7\equiv1\pmod\pi$, give
\[
(I+\rho+\rho^{-1})L_n\subseteq\pi\RR^{\Omega_n}.
\]
Multiplication by $\rho$ and saturation yield
\[
(\rho^2+\rho+I)L_n\subseteq\pi L_n.
\]
Subtracting $3\rho$ gives $(\rho-I)^2L_n\subseteq\pi L_n$.
Lemma~\ref{lem:local-torsion} applies because $\rho^n=I$. Extending
scalars gives $\rho^3=I$ on $W_n$. If $3\nmid n$, its order divides
both three and $n$, proving the last assertion.
\end{proof}

\subsection{The trace of rotation}\label{sec:trace}
For this calculation let $H$ have arbitrary order $q$. Associate to
each pair $a,b$ a unit vector $\xi_{ab}$ and a rank-one projection
$Q_{ab}$ by
\[
\xi_{ab}=q^{-1/2}(\overline{H_{ay}}H_{by})_{y\in X},\qquad
Q_{ab}=|\xi_{ab}\rangle\langle\xi_{ab}|.
\]
These row-quotient projections are those of Banica and Nicoara
\cite[Definition~2.3]{BN07}, with the quotient $H_b/H_a$. Fixing either index gives an orthonormal basis. The block matrix
$\mathsf V=(Q_{ab})\in M_q(\C)\otimes M_q(\C)$ is therefore unitary,
since
\[
\sum_bQ_{ab}Q_{cb}=\delta_{ac}I,\qquad
\sum_aQ_{ab}Q_{ad}=\delta_{bd}I.
\]
In terms of these vectors the profile is given by
\[
\langle\xi_{ab},\xi_{cd}\rangle=q^{-1}S(a,b,c,d).
\]

\begin{proposition}\label{prop:trace}
For a Hadamard subfactor of index $q$ and every $n\geq1$ one has
$|\Tr(\rho|P_n)|\leq q$. In particular, for Petrescu's matrix,
\begin{equation}\label{eq:trace-bound}
\left|\Tr(\rho|P_n)\right|\leq7.
\end{equation}
\end{proposition}
\begin{proof}
Let $\tau_q=q^{-1}\Tr_q$ on $M_q(\C)$. The trace of a cyclic product
of rank-one projections is the product of the successive inner
products. By \eqref{eq:Tn} this gives
\begin{equation}\label{eq:moment}
T_n(\mathbf a,\mathbf b)
=\tau_q(Q_{a_1b_1}\cdots Q_{a_nb_n}).
\end{equation}
In $M_q(\C)^{\otimes(r+1)}$, number the factors from $0$ to $r$ and
set
\[
\mathsf V^{[r]}=\mathsf V_{0,1}\mathsf V_{0,2}\cdots\mathsf V_{0,r}.
\]
Here $\mathsf V_{0,\ell}$ acts on factors $0,\ell$. This product
is unitary. Its matrix entries in factor zero are
\[
V^{[r]}_{ab}=\sum_{c_1,\ldots,c_{r-1}}
Q_{ac_1}\otimes Q_{c_1c_2}\otimes\cdots\otimes Q_{c_{r-1}b},
\]
with $V^{[1]}_{ab}=Q_{ab}$. Put $\tau_r=\tau_q^{\otimes r}$.
For the powers of $T_n$ we use Banica's moment formula
\cite[Proposition~3.2]{B14}, with the row-quotient convention above. In the present coordinates it follows by
expanding \eqref{eq:moment} and summing over the $r-1$ intermediate words,
\[
(T_n^r)(\mathbf a,\mathbf b)
=\tau_r(V^{[r]}_{a_1b_1}\cdots V^{[r]}_{a_nb_n}).
\]
Closing the outer indices cyclically gives the trace
\begin{align*}
\Tr(\rho T_n^r)
&=\sum_{a_1,\ldots,a_n}
\tau_r(V^{[r]}_{a_na_1}V^{[r]}_{a_1a_2}\cdots
V^{[r]}_{a_{n-1}a_n})\\
&=(\Tr_q\otimes\tau_r)((\mathsf V^{[r]})^n).
\end{align*}
The last expression is $q$ times a normalized trace of a unitary,
so its absolute value is at most $q$. By
Proposition~\ref{prop:transfer}, $T_n^r$ converges to the orthogonal
projection onto $\mathcal F_n$. The cyclic product \eqref{eq:Tn}
commutes with $\rho$. Passing to the limit gives the assertion.
\end{proof}

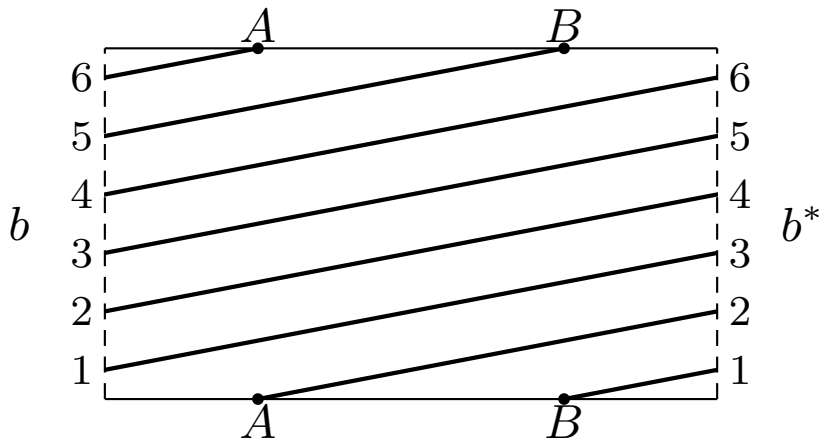
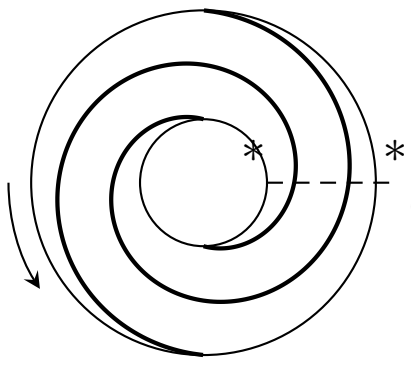
\begin{figure}[htbp]
\centering
\begin{tikzpicture}[font=\footnotesize,>=stealth]
\begin{scope}[shift={(-3.0,0)}]
 \draw (0,0) circle (1.14);
 \draw (0,0) circle (.42);
 \foreach \offset in {90,270}{
   \draw[thick,domain=0:1,samples=110,smooth,variable=\t]
      plot ({(.42+.72*\t)*cos(\offset+360*\t)},
            {(.42+.72*\t)*sin(\offset+360*\t)});
 }
 \draw[densely dashed] (.42,0)--(1.25,0);
 \node[above] at (.33,0) {$*$};
 \node[above right] at (1.05,0) {$*$};
 \node[right] at (1.23,-.13) {cut};
 \draw[->] (-1.29,0) arc[start angle=180,end angle=213,radius=1.29];
 \node at (0,-1.62) {The winding annulus $u$};
\end{scope}
\begin{scope}[shift={(.25,-1.16)},x=4.05cm,y=2.32cm]
 \begin{scope}
  \clip (0,0) rectangle (1,1);
  \foreach \offset in {.25,.75}{
   \foreach \k in {0,1,2,3}{
    \draw[thick] ({\offset-\k},0)--({3+\offset-\k},1);
   }
  }
 \end{scope}
 \draw (0,0)--(1,0) (0,1)--(1,1);
 \draw[densely dashed] (0,0)--(0,1) (1,0)--(1,1);
 \foreach \j in {1,...,6}{
   \pgfmathsetmacro{\h}{(2*\j-1)/12}
   \node[anchor=east,inner sep=2pt] at (0,\h) {\scriptsize$\j$};
   \node[anchor=west,inner sep=2pt] at (1,\h) {\scriptsize$\j$};
 }
 \foreach \ab/\letter in {.25/A,.75/B}{
  \fill (\ab,0) circle (1.1pt);
  \fill (\ab,1) circle (1.1pt);
  \node[below,inner sep=1pt] at (\ab,0) {$\letter$};
  \node[above,inner sep=1pt] at (\ab,1) {$\letter$};
 }
 \node at (-.14,.5) {$b$};
 \node at (1.14,.5) {$b^*$};
 \node at (.5,-.21) {The cut annulus $u^3$};
\end{scope}
\end{tikzpicture}
\caption{The stars mark the intervals used for the radial cut.
In the right panel the banks of the cut are vertical. Closing the
two-box joins $A$ to $A$ and $B$ to $B$. A strand from the left
bank at~$1$ reaches the right bank at~$3$. The strand from~$5$
continues through the closure to~$1$.}
\label{fig:affine-cut}
\end{figure}

\begin{remark}\label{rem:affine-trace}
The trace also admits an affine planar algebra interpretation. Let
$\delta$ be the modulus of the spherical subfactor planar algebra $P$.
Write $\operatorname{tr}_k$ for its normalized picture trace on $P_{k,+}$.
Das, Ghosh and Gupta construct a conditional expectation
$\gamma_{+1}:AP_{+1,+1}\to P_{2,+}$ and a tracial state
$\omega_{+1}=\operatorname{tr}_2\circ\gamma_{+1}$, also at infinite depth
\cite[Lemma~3.3 and Corollary~3.9]{DGG14}.

Let $u$ be the affine annulus on two through strings, each making one
full turn, with the direction chosen to agree with $\rho$.
Choose the radial cut through the two distinguished intervals.
In $u^n$ each string crosses the cut $n$ times. Number the $2n$
crossings on each bank from the inner circle to the outer circle.
The two strings alternate on the cut. Hence a string leaving crossing
$j$ reaches crossing $j+2$ on the other bank. For the last two
crossings this numbering is taken modulo $2n$, since the picture trace
joins the two outer endpoints to the corresponding inner endpoints.
Thus the tangle between the two banks after closure is precisely
two-click rotation (Figure~\ref{fig:affine-cut}).

Take an orthonormal basis $\mathcal B_n$ of $P_{n,+}$ for
$\langle b,c\rangle=\operatorname{tr}_n(b^*c)$.
In the formula for $\gamma_{+1}$, the cut with $2n$ intersections
inserts $b$ and $b^*$ on its banks, sums over $b\in\mathcal B_n$,
and contributes the factor $\delta^{-n}$. Closing the output
two-box for $\operatorname{tr}_2$ contributes $\delta^{-2}$.
For fixed $b$ the unnormalized closed picture is the closure of
$b^*\rho(b)$ and has value $\delta^n\operatorname{tr}_n(b^*\rho(b))$.
Consequently
\[
 \omega_{+1}(u^n)
 =\delta^{-2}\delta^{-n}\sum_{b\in\mathcal B_n}
       \delta^n\operatorname{tr}_n(b^*\rho(b))
 =\delta^{-2}\Tr(\rho|P_n).
\]
Reflection reverses the winding, so $u^*=u^{-1}$.
For a Hadamard subfactor $\delta^2=q$, and
$|\omega_{+1}(u^n)|\leq1$ recovers Proposition~\ref{prop:trace}.
\end{remark}

\section{The standard invariant}\label{sec:annular}
\subsection{Annular multiplicities}
Let
\[
b_n=\dim P_n\quad(n\geq0),\qquad
a_n=\dim W_n\quad(n\geq1).
\]
For modulus $\delta>2$, the irreducible Hilbert annular
Temperley--Lieb modules of lowest weight $k>0$ are indexed by
the rotation eigenvalues $\omega$, $\omega^k=1$. The lowest-weight
space of $V^{k,\omega}$ is one-dimensional and its higher dimensions
are
\[
\dim V_n^{k,\omega}=\binom{2n}{n-k}\qquad(n\geq k).
\]
The common cap kernel is the space of lowest-weight vectors. Hence
$a_k$ counts the summands of lowest weight $k$, including
multiplicities, in the orthogonal annular decomposition
\cite[Remark~3.6, Corollary~3.9, Lemma~5.1, Definition~5.2,
Proposition~5.3, Theorem~5.7 and Corollary~5.8]{J01}.
The module generated by the empty diagrams is Temperley--Lieb
\cite[Proposition~5.13]{J01}. It exhausts weight zero since the two
zero-box spaces are one-dimensional. With $a_1=0$ this gives
\begin{equation}\label{eq:annular-dimension}
b_n=\Cat_n+\sum_{k=2}^n\binom{2n}{n-k}a_k.
\end{equation}
This is the coefficient identity of \cite[Theorem~6.3]{J01}.

\begin{lemma}\label{lem:odd-trace}
For odd $n\geq3$ one has
\begin{equation}\label{eq:odd-trace}
\Tr(\rho|P_n)=2+\Tr(\rho|W_n).
\end{equation}
For the Petrescu subfactor it follows that
\begin{equation}\label{eq:bounded-subsequence}
0\leq a_n\leq5\qquad(n\geq3\text{ odd},\ 3\nmid n).
\end{equation}
\end{lemma}
\begin{proof}
In a summand of lowest weight $0<k<n$, use the annular diagram
basis of \cite[Definition~5.2]{J01}. For a summand with rotation eigenvalue $\omega$, choose one diagram
from each orbit of inner rotation and project it onto the
$\omega$-eigenspace. The resulting vectors form a basis. All diagrams
in a given orbit have the same set of $2k$ outer through-string
endpoints, and two-click rotation of the outer boundary permutes
these basis vectors up to scalars. A basis vector contributing to the
trace must therefore have its set of outer through-string endpoints
invariant under two-click rotation. The latter has two orbits of
size $n$. Since $n$ is odd, a union of these orbits cannot have
even cardinality strictly between zero and $2n$. Thus all these
summands have trace zero.

A Temperley--Lieb basis diagram has an adjacent pair. If rotation
fixes the diagram, the orbit of this pair gives all its pairs.
For $n\geq3$ there are just two such diagrams, according to the
parity of the first endpoint of an adjacent pair. The summands of
weight $n$ contribute precisely $\Tr(\rho|W_n)$. This proves
\eqref{eq:odd-trace}. If $3\nmid n$, Theorem~\ref{thm:rotation}
and Proposition~\ref{prop:trace} give $2+a_n\leq7$.
\end{proof}

\subsection{The norm of the principal graph}
The next assertion concerns an arbitrary irreducible subfactor.
We use $a_n,b_n$ in the sense of \eqref{eq:annular-dimension}.

\begin{proposition}\label{prop:norm}
Let $P$ be the planar algebra of an irreducible finite-index subfactor
of modulus $\delta>2$. If the annular multiplicities $a_n$ are
bounded on an unbounded sequence of degrees, the principal graph
has norm at most two.
\end{proposition}
\begin{proof}
Write $A$ for the adjacency operator of the principal graph $\Gamma$
and $*$ for its root. The positive dimension function satisfies
$\sum_w A_{vw}d_w=\delta d_v$, so the weighted Schur test gives
$\|A\|\leq\delta$. Let $\mu$ be the spectral measure of $A^2$ at
$e_*$. This is a probability measure on $[0,\delta^2]$. Counting
closed paths in the tower \cite[Section~2.6(iii)]{J12} gives
\[
 b_j=\langle A^{2j}e_*,e_*\rangle
     =\int x^j\,d\mu(x).
\]

To recover $a_n$ from these moments, put $q_0=2$, $q_1=x-2$ and
\[
 q_{n+1}=(x-2)q_n-q_{n-1}.
\]
Then $q_n(2+2\cos\theta)=2\cos(n\theta)$. The expansion at
$x=2+z+z^{-1}$ gives
\begin{equation}\label{eq:binomial}
 x^n=\binom{2n}{n}+\sum_{k=1}^n\binom{2n}{n-k}q_k(x).
\end{equation}
Now $\int q_1\,d\mu=-1$ since $b_0=b_1=1$. Integrating
\eqref{eq:binomial}, the first two coefficients contribute
$\binom{2n}{n}-\binom{2n}{n-1}=\Cat_n$. The remaining triangular
system has diagonal entries one. Comparing it with
\eqref{eq:annular-dimension} gives
\begin{equation}\label{eq:qmoment}
 a_n=\int q_n(x)\,d\mu(x)\qquad(n\geq2).
\end{equation}

On $[0,4]$ one has $q_n\geq-2$, whereas for $x>4$,
\[
 q_n(x)=\lambda(x)^n+\lambda(x)^{-n},\qquad
 \lambda(x)=\frac{x-2+\sqrt{x(x-4)}}2>1.
\]
If $\mu((4,\delta^2])>0$, choose $c>4$ such that
$\mu([c,\delta^2])=\varepsilon>0$. Since $\lambda$ is increasing,
\eqref{eq:qmoment} gives
\[
 a_n\geq\varepsilon\lambda(c)^n-2.
\]
This would force $a_n$ to tend to infinity. Thus $\mu$ is supported
in $[0,4]$ and
\[
 b_j=\int x^j\,d\mu(x)\leq4^j\qquad(j\geq1).
\]
To pass from the root to the whole graph, let $G$ be a finite
connected subgraph containing $*$, with positive unit Perron vector
$v$. The spectral decomposition of $A_G$ and inclusion of the
rooted walks in those of $\Gamma$ give
\[
 |v_*|^2\|A_G\|^{2j}\le(A_G^{2j})_{*,*}\le b_j\le4^j.
\]
Taking $2j$th roots gives $\|A_G\|\leq2$. These finite subgraphs
exhaust $\Gamma$, proving $\|A\|\leq2$.
\end{proof}

\subsection{The principal and dual graphs}
\begin{lemma}\label{lem:graph}
Suppose an irreducible subfactor has modulus $\delta>2$, principal
graph of norm at most two, and $\dim P_2=2$. Its rooted principal
graph is $A_\infty$.
\end{lemma}
\begin{proof}
On a finite graph the positive dimension vector would be an
eigenvector of eigenvalue $\delta>2$. The graph is therefore
infinite. The connected infinite bipartite graphs of norm at most
two are $A_\infty$, $A_{\infty,\infty}$ and $D_\infty$
\cite[Section~2.10]{P23}. Irreducibility gives a root of degree one.
This excludes $A_{\infty,\infty}$. From either leaf of $D_\infty$
there are three closed walks of length four, whereas $\dim P_2=2$.
The remaining graph is $A_\infty$ rooted at its endpoint.
\end{proof}

\begin{proof}[Proof of Theorem~\ref{thm:main}]
By Lemma~\ref{lem:odd-trace}, $a_{6m+1}\leq5$ for $m\geq1$.
Proposition~\ref{prop:norm} bounds the graph norm by two, and
Lemmas~\ref{lem:low} and~\ref{lem:graph} identify the rooted
principal graph. Its closed walks of length $2n$ are Dyck paths,
so $\dim P_n=\Cat_n$. The canonical inclusion
$\TL_n(\sqrt7)\subseteq P_n$ is consequently an equality.
One-click rotation carries the positive box space isomorphically
onto the negative box space and preserves the Temperley--Lieb
subspaces \cite[Section~2.6(iv)]{J12}. This gives equality in the
opposite shading, identifies the dual graph, and proves
\eqref{eq:dimensions} and \eqref{eq:dual-dimensions}.
\end{proof}

\begin{corollary}\label{cor:rank}
For every $n\geq1$,
\[
\rank(I-T_n)=7^n-\Cat_n.
\]
Moreover $W_n=0$.
\end{corollary}
\begin{proof}
Apply rank-nullity to \eqref{eq:fixed}. In \eqref{eq:annular-dimension}
the Catalan term accounts for all of $b_n$, so the positive-weight
multiplicities vanish. The case $n=1$ follows from
Proposition~\ref{prop:transfer}.
\end{proof}
In the indexing of the Jones tower the first dimensions are
\[
\begin{array}{c|rrrrrrrrr}
k&0&1&2&3&4&5&6&7&8\\\hline
\dim(N'\cap M_k)&1&2&5&14&42&132&429&1430&4862.
\end{array}
\]

\section{Specialization in Petrescu's family}\label{sec:family}
Let $\T=\{z\in\C:|z|=1\}$ and $w=e^{\pi i/3}$. We use the
family of \cite[Section~4.3]{N10} in the form
\begin{equation}\label{eq:family-matrix}
H(z)=\begin{pmatrix}
1&1&1&1&1&1&1\\
1&zw&zw^4&w^5&-1&-1&w\\
1&zw^4&zw&-1&w^5&-1&w\\
1&w^5&-1&z^{-1}w&z^{-1}w^4&w&-1\\
1&-1&w^5&z^{-1}w^4&z^{-1}w&w&-1\\
1&-1&-1&w&w&w^4&w^5\\
1&w&w&-1&-1&w^5&w^4
\end{pmatrix},\qquad z\in\T.
\end{equation}
The rows and columns are indexed by $X_7=\{0,\ldots,6\}$.
The identities $\overline z=z^{-1}$ and $w^2=w-1$ give
$H(z)H(z)^*=7I_7$. At $z=1$ one recovers $P_7$, which is
equivalent to \eqref{eq:H} by Lemma~\ref{lem:matrix}.

The Jones tower is
\[
N_z\subset M_z=M_{z,0}\subset M_{z,1}\subset\cdots.
\]
Write $P_n(z)=N_z'\cap M_{z,n-1}$ and $P_0(z)=\C$. The
Jones projections give the inclusion
\[
\TL_n(\sqrt7)\subseteq P_n(z)
\]
for all $z$. In particular $\dim P_n(z)\geq\Cat_n$. At the base
point Theorem~\ref{thm:main} gives
\begin{equation}\label{eq:family-base}
P_n(1)=\TL_n(\sqrt7)\qquad(n\geq0).
\end{equation}

\begin{theorem}\label{thm:family}
For each $n\geq1$ the exceptional set
\[
\mathcal E_n=\{z\in\T:P_n(z)\ne\TL_n(\sqrt7)\}
\]
is finite and consists of algebraic numbers. The union
$\mathcal E=\bigcup_{n\geq1}\mathcal E_n$ is therefore at most
countable and excludes $1$. For $z\in\T\setminus\mathcal E$, in
particular for every transcendental phase, the standard invariant
is $\TL(\sqrt7)$ and both principal graphs are $A_\infty$.
\end{theorem}

Put $K=\Q(w)$ and give $K[t,t^{-1}]$ the involution which sends
$h(t)=\sum_jc_jt^j$ to
\[
h^{\#}(t)=\sum_j\overline{c_j}\,t^{-j}.
\]
Evaluation at any $z\in\T$ intertwines this involution with complex
conjugation.

\begin{lemma}\label{lem:family-equations}
For $n\geq1$ there is a matrix
\[
\mathsf A_n(t)\in M_{7^n}\bigl(K[t,t^{-1}]\bigr)
\]
such that $\ker\mathsf A_n(z)\cong P_n(z)$ for every $z\in\T$.
\end{lemma}
\begin{proof}
Denote the entries of \eqref{eq:family-matrix} by $h_{ay}(t)$ and put
\[
S_t(a,b,c,d)=\sum_{y\in X_7}
h_{ay}(t)h_{by}^{\#}(t)h_{cy}^{\#}(t)h_{dy}(t).
\]
For words $\mathbf a,\mathbf b\in X_7^n$ set
\begin{equation}\label{eq:family-transfer}
T_n(t)_{\mathbf a,\mathbf b}
=7^{-n-1}\prod_{j=1}^n S_t(a_j,b_j,a_{j+1},b_{j+1}),
\qquad \mathsf A_n(t)=I-T_n(t),
\end{equation}
with subscripts read cyclically. These matrices are Laurent
polynomial over $K$. On $\T$ their coefficients specialize to
the profiles and transfer matrices of $H(z)$, so
Proposition~\ref{prop:transfer} identifies the asserted kernel.
\end{proof}

\begin{proof}[Proof of Theorem~\ref{thm:family}]
Let $r_n=7^n-\Cat_n$. At $t=1$ the preceding lemma and
\eqref{eq:family-base} give rank $r_n$. Choose an $r_n$-rowed minor
$\Delta_n(t)$ nonzero there. Wherever $\Delta_n(z)\ne0$ we have
\[
\dim P_n(z)=7^n-\rank\mathsf A_n(z)
\leq7^n-r_n=\Cat_n.
\]
The reverse inequality holds throughout the family. Hence
\[
\mathcal E_n\subseteq\{z\in\T:\Delta_n(z)=0\}.
\]
Multiplication by a power of $t$ turns $\Delta_n$ into a nonzero
polynomial over $K$ with the same zeros on $\T$. This proves the
assertion at fixed $n$. Taking the union over $n$ gives the stated
exceptional set, and \eqref{eq:family-base} excludes $1$. For every
remaining parameter the positive box spaces are Temperley--Lieb.
Rotation gives the opposite shading as in Theorem~\ref{thm:main}.
\end{proof}

\begin{remark}\label{rem:family-local}
At any prescribed depth $m$ there is an arc about $1$ on which
$P_n(z)=\TL_n(\sqrt7)$ for every $n\leq m$. One need only avoid
the finite union $\bigcup_{n\leq m}\mathcal E_n$. The arc may
depend on $m$. Finiteness or emptiness of the full exceptional
set does not follow from this argument.
\end{remark}

\section{Algebraic parameters}\label{sec:algebraic-family}
We retain the family \eqref{eq:family-matrix}. For an algebraic phase
the reduction of its entries can be used in place of specialization.

\begin{theorem}\label{thm:algebraic-family}
Let $z\in\overline\Q\cap\T$, $F=\Q(w,z)$, and let $\OO_F$ be its
ring of integers. Suppose that at some prime $\mathfrak p$ over three,
\begin{equation}\label{eq:algebraic-criterion}
z\in R=(\OO_F)_{\mathfrak p},\qquad
z\equiv1\ \text{or}\ -1\pmod{\mathfrak pR}.
\end{equation}
Then
\[
P_n(z)=\TL_n(\sqrt7)\qquad(n\geq0).
\]
The full standard invariant is $\TL(\sqrt7)$ and both principal
graphs are $A_\infty$.
\end{theorem}

The condition makes $z$ a local unit. The proof allows arbitrary
ramification at $\mathfrak p$ and uses conjugation only through
the identity $\overline z=z^{-1}$.

\begin{lemma}\label{lem:prime-to-three}
Let $(R,\mathfrak m)$ be a discrete valuation ring of residue
characteristic three, $L$ a finite free $R$-module, and
$B\in\operatorname{Aut}_R(L)$. If
\[
B^n=I,\qquad 3\nmid n,\qquad (B-I)^2L\subseteq\mathfrak mL,
\]
then $B=I$.
\end{lemma}
\begin{proof}
The minimal polynomial of the reduction of $B$ divides both
$(t-1)^2$ and $t^n-1$. The latter is separable in the residue
field, so $B$ reduces to $I$. Since $n$ is a unit in $R$,
\[
E=\frac1n\sum_{j=0}^{n-1}B^j
\]
is an idempotent with $BE=E$. It reduces to $I$ and hence has
unit determinant. Thus $E=I$ and $B=I$.
\end{proof}

\begin{proof}[Proof of Theorem~\ref{thm:algebraic-family}]
Put $\mathfrak m=\mathfrak pR$. Apply the permutations and signs
of Lemma~\ref{lem:matrix} to obtain $\widehat H(z)$, with indices
in $X=\{\infty\}\sqcup\Z/6\Z$. In terms of the coordinate vectors
$\eta_a$ let
\[
v=\eta_0-\eta_2,\qquad u=\eta_3-\eta_5.
\]
Comparison of the two variable blocks in \eqref{eq:family-matrix}
with their values at $z=1$ gives
\begin{equation}\label{eq:family-rank-two}
\widehat H(z)=H+w(z-1)vv^{\mathsf T}
                 +w(z^{-1}-1)uu^{\mathsf T},
\end{equation}
where $H$ is the sixth-root matrix \eqref{eq:H}.

We have $w\equiv-1$ since $(1+w)^2=3w$. If $z\equiv1$, reduction
of \eqref{eq:family-rank-two} gives $I_7+J_7$. If $z\equiv-1$,
let
\[
V=I_7-vv^{\mathsf T}-uu^{\mathsf T}.
\]
This permutation transposes $0,2$ and $3,5$, so $V^2=I_7$ and
$VJ_7=J_7$. In this case $\widehat H(z)\equiv V+J_7$ and
$V\widehat H(z)\equiv I_7+J_7$. Choose the resulting equivalent
representative $\widetilde H$ in the two cases. Each entry is a unit
with residue $\pm1$. Its conjugate is its inverse and has the same
residue. We have proved
\begin{equation}\label{eq:family-local-reduction}
\widetilde H\equiv\overline{\widetilde H}
\equiv I_7+J_7\pmod{\mathfrak m}.
\end{equation}

By Lemma~\ref{lem:profile-congruence}, its profile satisfies
\begin{equation}\label{eq:family-profile-reduction}
\widetilde S(a,b,c,d)\equiv
-1+\delta_{ab}+\delta_{ad}+\delta_{cb}+\delta_{cd}
\pmod{\mathfrak m}
\end{equation}
whenever $a\ne c$ and $b\ne d$.
Let $\widetilde K_n$ be the cyclic profile product for
$\widetilde H$ on $\Omega_n$. The cap equations in
Lemma~\ref{lem:cyclic-expansion} then give
\begin{equation}\label{eq:family-cyclic-congruence}
\widetilde K_nf\equiv(2I+\rho+\rho^{-1})f
\pmod{\mathfrak m}
\end{equation}
for $R$-valued $f\in\mathcal U_n$. This computation uses only
\eqref{eq:family-profile-reduction} in the residue field.

For $n\geq2$ let $W_{n,F}(z)$ denote the solutions over $F$ of the
flatness and cap equations. Put
\[
L_n=W_{n,F}(z)\cap R^{\Omega_n}.
\]
As in Theorem~\ref{thm:rotation}, this is a saturated free
$R$-lattice in the solution space. It is stable under $\rho$ and
has complexification $W_n(z)$. On $L_n$ the flatness equation
is $\widetilde K_nf=7^{n+1}f$. Since $7\equiv1$, substitution in
\eqref{eq:family-cyclic-congruence}, multiplication by $\rho$, and
saturation give
\[
(\rho^2+\rho+I)L_n\subseteq\mathfrak mL_n.
\]
Subtracting $3\rho$ gives the hypothesis of
Lemma~\ref{lem:prime-to-three} for $B=\rho$ in degrees prime to
three. Thus
\begin{equation}\label{eq:family-rotation}
\rho|_{W_n(z)}=I\qquad(3\nmid n).
\end{equation}
The same identity holds for $n=1$ since $W_1(z)=0$.

For odd $n\geq3$ in these degrees, Proposition~\ref{prop:trace}
and \eqref{eq:odd-trace} give
\[
\Tr(\rho|_{P_n(z)})=2+\dim W_n(z),\qquad
\dim W_n(z)\leq5.
\]
Proposition~\ref{prop:norm} bounds the graph norm by two.
Equation~\eqref{eq:family-profile-reduction} gives nonzero entries
between disjoint pairs in the off-diagonal block of $T_2$. As in
Lemma~\ref{lem:low} this block is irreducible, so $\dim P_2(z)=2$.
Since $T_1=J_7/7$, the inclusion is irreducible.
Lemma~\ref{lem:graph} identifies its principal graph as $A_\infty$.
The Catalan dimensions and one-click rotation now give equality
with Temperley--Lieb in both shadings.
\end{proof}

\begin{corollary}\label{cor:algebraic-parameters}
The standard invariant is $\TL(\sqrt7)$, with both principal graphs
$A_\infty$, for every $z\in\Q(w)\cap\T$ and every root of unity of
order dividing $2\cdot3^r$ for some $r\geq0$.
\end{corollary}
\begin{proof}
In $\Z[w]$ the unique prime $(1+w)$ above three is fixed by
conjugation. For $z\in\Q(w)\cap\T$ its valuation consequently gives
\[
0=\nu(z)+\nu(\overline z)=2\nu(z).
\]
Thus $z$ is a unit with nonzero residue in $\F_3$, hence residue
$\pm1$. For $z^{2\cdot3^r}=1$, the phase is a unit at every
finite prime and its reduction at a prime over three satisfies
\[
0=z^{2\cdot3^r}-1=(z^2-1)^{3^r}
\]
in the residue field. Again $z$ reduces to $\pm1$, so the theorem
applies in both cases.
\end{proof}

For example, $(3+5w)/7$ is a phase in $\Q(w)$ of infinite order.
Indeed its sum with its conjugate is $11/7$, which is not an
algebraic integer. It lies outside $\mathcal E$, as do all roots
of unity of orders $3^r$ and $2\cdot3^r$.

\section{Extensions to other Hadamard matrices}\label{sec:extensions}
Let $H\in M_q(\C)$ be a complex Hadamard matrix, so that
$|H_{ay}|=1$ and $HH^*=qI_q$. Write $\mathrm{BH}(q,k)$ for those
matrices whose entries are $k$th roots of unity. We keep the
notation $P_n=N'\cap M_{n-1}$. To prove that the standard invariant
is $\TL(\sqrt q)$ it suffices, for $q>4$, to establish
$\dim P_2=2$ and boundedness of $\dim W_n$ along an unbounded
sequence of degrees. Indeed $P_1=\C$, and Proposition~\ref{prop:norm}
and Lemma~\ref{lem:graph} then apply as in Theorem~\ref{thm:main}.
In all such cases $\dim(N'\cap M_j)=\Cat_{j+1}$,
$\dim(M'\cap M_j)=\Cat_j$, and both principal graphs are $A_\infty$.

We may pass to any equivalent Hadamard matrix. In \eqref{eq:Cn} row
phases cancel around the cycle, column phases cancel in pairs, and
permutations relabel the words. Conjugation conjugates $C_n$.
Transposition of $H$ transposes $C_n$ after reversal of the two cyclic
words. Thus conjugation, transposition and adjoint also preserve its
singular values and the dimensions of the fixed spaces of $C_nC_n^*$.

\subsection{Arithmetic conditions on the entries}
\begin{theorem}\label{thm:ext-classes}
Let $H$ be a complex Hadamard matrix of order $q$. Under any of the
following hypotheses its associated subfactor has standard invariant
$\TL(\sqrt q)$.
\begin{enumerate}
\item $H$ is real, $q>4$ and $q\equiv4\pmod8$.
\item $H\in M_q(\Q(i))$, $q>2$ and $q\equiv2\pmod4$. This
includes all of $\mathrm{BH}(q,4)$ in these orders.
\item $H\in\mathrm{BH}(q,3)$, $q>3$, $3\mid q$ and $9\nmid q$.
\item $H$ has an algebraic representative with local unit entries
reducing to $J_q-2I_q$ at a prime of odd residue characteristic $p$,
where $q>4$ and $q\equiv4\pmod p$.
\end{enumerate}
\end{theorem}

For a prime power $s\equiv1\pmod4$, the Paley type~II Hadamard
matrix is real of order $q=2(s+1)\equiv4\pmod8$. Since $s\geq5$,
part~(1) gives standard invariant $\TL(\sqrt{2(s+1)})$. This proves
Montgomery's Conjecture~3.5 \cite{M21}.

Petrescu's reduction $I_7+J_7$ in characteristic three is an instance
of (4). For an example of (2), let $C$ be any real symmetric
conference matrix of order $q\equiv2\pmod4$, $q>2$. The diagonal
is zero, the off-diagonal entries are signs and $C^2=(q-1)I_q$.
Hence $C+iI_q$ has fourth-root entries and
$(C+iI_q)(C-iI_q)=qI_q$. In particular the theorem applies to
the Paley conference construction in these orders.

For $a\ne c$, $b\ne d$, put
\[
 D=\delta_{ab}+\delta_{cd}+\delta_{ad}+\delta_{bc},\qquad
 E=\delta_{ab}+\delta_{cd}-\delta_{ad}-\delta_{bc}.
\]
We call the profile on such pairs the restricted profile. In each
case the argument depends on the reduction of $S/\alpha$ in the
local ring of the indicated prime, for the $\alpha$ of
Table~\ref{tab:ext-profiles}. Here
$S(a,b,c,d)=\sum_yH_{ay}H_{by}^{-1}H_{cy}^{-1}H_{dy}$.
\begin{table}[htbp]
\caption{The restricted profiles. The entry in the last column is taken
in the residue field of characteristic $p$.}\label{tab:ext-profiles}
\centering
\begin{tabular}{@{}llccc@{}}\toprule
Entries or reduction & Order & $p$ & $\alpha$ & $S/\alpha$ \\ \midrule
Real & $q=4m$, $m$ odd & 2 & 4 & $1+D$\\
Gaussian rational & $q=2m$, $m$ odd & 2 & 2 & $1+D$\\
Third roots & $q=3m$, $3\nmid m$ & 3 & 3 & $-m(1-E)$\\
$J_q-2I_q$ & $q\equiv4\pmod p$ & $p$ odd & 1 & $4(-1+D)$\\
\bottomrule\end{tabular}
\end{table}

\subsection{Congruences for the restricted profile}
\emph{Real matrices.}
Montgomery obtained the following congruence for four distinct rows
and used it to prove $\dim P_2=2$ \cite[Proposition~3.2.9]{M22}.
We record the calculation for use in the cyclic product at degree $n$.
For four distinct indices multiply each column by its entry in
row $a$. This does not change $S$. Write the other three rows as
$1-2x$, $1-2y$, $1-2z$. Their supports of negative entries have
size $q/2$, with pairwise intersections of size $q/4$, by row
orthogonality. Expanding the product gives
$S=q-8\sum_jx_jy_jz_j\equiv q\pmod8$. For pairs sharing precisely
one label, $S$ is an inner product of distinct rows and vanishes.
For equal unordered pairs $S=q$. Dividing by four gives
$S/4\equiv1+D\pmod2$ in all three cases.

\emph{Gaussian-rational matrices.}
Use a dephased representative and put $\pi=1-i$,
$R=\Z[i]_{(\pi)}$. Every phase in $\Q(i)$ is a local unit congruent
to a fourth root modulo $4R$. To see this, write the phase as
$(a+bi)/c$ with primitive integer coordinates. The equation
$a^2+b^2=c^2$ modulo four gives $c$ odd and $a,b$ of opposite
parity. Modulo eight, the even one is divisible by four, while
the odd one divided by $c$ is $\pm1$ modulo four.

Choose fourth roots $G_{ay}=i^{e_{ay}}\equiv H_{ay}\pmod{4R}$,
$e_{ay}\in\{0,1,2,3\}$, with first exponent row zero. Inversion
preserves this congruence, so $\sum_yi^{e_{ay}-e_{by}}\in4R$
for distinct rows. Comparing with the first row, the imaginary part
shows that each row has an even number of odd exponents. If $r$ is
an integer vector of length $2m$ with even sum, expand
$(1-\pi)^{r_y}$ through degree two and use $\pi^2=-2i$. One obtains
\[
 \frac12\sum_y i^{r_y}\in R,\qquad
 \frac12\sum_y i^{r_y}\equiv m+\sum_y\binom{r_y}{2}\pmod\pi.
\]
Let $v_a=e_a\bmod2$ and $Q_a=\sum_y\binom{e_{ay}}2\bmod2$.
Applied to $r=e_a-e_b$, the identity gives, for all $a,b$,
\[
 v_a\mathbin{\cdot}v_b=m(1+\delta_{ab})+Q_a+Q_b\quad\text{in }\F_2.
\]
Here $\binom{-t}{2}=\binom t2+t$ and the row sums are even.
For $r=e_a-e_b-e_c+e_d$, the binomial sum contains four terms
$Q$ and six pairwise inner products, counting the four positions
with multiplicity. Substitution of the last formula cancels the
$Q$'s. The remaining terms are $mD$ since $a\ne c$, $b\ne d$.
It follows that $S_G/2\equiv m(1+D)=1+D\pmod\pi$. The original
profiles agree with those of $G$ modulo $4R$, which proves the
second reduction.

\emph{Cube roots.}
Write $q=3m$, $\omega=e^{2\pi i/3}$, $H_{ay}=\omega^{e_{ay}}$,
with first exponent row zero. Since $1+\omega+\omega^2=0$ is the
minimal relation, orthogonality makes every difference of distinct
exponent rows take each residue $m$ times. In $\F_3$ such a
difference has sum zero and squared length $2m$. Polarization gives
diagonal entries $2m$ and off-diagonal entries $m$ in the Gram matrix
of the nonzero rows. Including the zero row, this gives for
$\sum_j t_j=0$
\[
 \left(\sum_j t_je_j\right)\mathbin{\cdot}
 \left(\sum_j t_je_j\right)=m\sum_j t_j^2.
\]
Let $t$ be the coefficients of $e_a-e_b-e_c+e_d$, and $n_j$ the
number of coordinates of this vector equal to $j\in\F_3$.
Its coordinate sum is zero and its length is divisible by three,
so $n_0\equiv n_1\equiv n_2\pmod3$. Thus
\[
 S=(n_0-n_2)+(n_1-n_2)\omega\in3\Z[\omega],\qquad
 S/3\equiv m-n_2\pmod{1-\omega}.
\]
Summing the indicator $(x^2-x)/2$ of $x=2$ gives
$n_2=2m\sum_jt_j^2$ in $\F_3$. On restricted pairs
$\sum_jt_j^2=4-2E$, which gives $S/3\equiv-m(1-E)$.

\emph{Reduction to $J_q-2I_q$.}
The reduced entries are signs. In the product of the four rows,
a column contributes $-1$ precisely on the symmetric difference
of $\{a,c\}$ and $\{b,d\}$. This set has $4-2D$ elements.
Hence the reduced sum is $q-8+4D=4(-1+D)$.

\begin{proof}[Proof of Theorem~\ref{thm:ext-classes}]
Let $(R,\mathfrak m)$ be the local ring and put $B_n=\alpha^{-n}K_n$
on $\Omega_n$. Its entries are integral by the calculation above.
For a flat lowest-weight vector the eigenvalue equation is
\begin{equation}\label{eq:ext-flat}
 B_nf=\frac{q^{n+1}}{\alpha^n}f\qquad(f\in W_n).
\end{equation}
As in Section~\ref{sec:algebraic-family}, the solutions over the
field of entries intersect $R^{\Omega_n}$ in a saturated,
rotation-stable lattice $L_n$ whose complexification is $W_n$.

Fix an odd $n\geq3$. In the expansion of $B_nf$ for $f\in L_n$,
Lemma~\ref{lem:cyclic-expansion} leaves four terms modulo
$\mathfrak m$. The constant choices of $\delta_{ab}$ and
$\delta_{cd}$ act as the identity, and the choices of $\delta_{ad}$
and $\delta_{bc}$ are the two shifts. The first two rows of the
table therefore give $2I+\rho+\rho^{-1}$. For cube roots the
two identity terms acquire the sign $(-1)^n$, giving
$(-m)^n(-2I+\rho+\rho^{-1})$. For the last row the result is
$4^n(2I+\rho+\rho^{-1})$.

The scalar in \eqref{eq:ext-flat} is $\alpha m^{n+1}$ in the first
three cases, hence zero in the residue field. In the fourth it is
$4^{n+1}$ in the residue field. Multiplying the resulting congruence
by $\rho$ and using saturation gives in all cases
\[
 (\rho-I)^2 L_n\subseteq\mathfrak mL_n.
\]
For $p\nmid n$ the proof of Lemma~\ref{lem:prime-to-three} applies
in residue characteristic $p$. Indeed $t^n-1$ is separable, so the
reduction of $\rho$ is $I$, and its averaging idempotent is invertible
over $R$. It follows that $\rho$ fixes $W_n$. Proposition~\ref{prop:trace}
and \eqref{eq:odd-trace} give
\[
 \dim W_n\leq q-2\qquad(n\geq3\text{ odd},\ p\nmid n).
\]
Thus Proposition~\ref{prop:norm} applies. For disjoint pairs all four
residues in Table~\ref{tab:ext-profiles} are nonzero. The stochastic
block of $T_2$ on distinct pairs has connected support, since
$q\geq6$ and any two pairs admit a third disjoint from both. This
block and the diagonal block each contribute one fixed vector, as
in Lemma~\ref{lem:low}. We have obtained $\dim P_2=2$ and the norm
bound, proving the assertion in both shadings.
\end{proof}

\subsection{Specialization and algebraic phases}\label{ssec:ext-affine}
\begin{proposition}\label{prop:ext-specialize}
Let $H(\mathbf z)$ be a Hadamard family of order $q>4$ on $\T^d$.
Suppose its entries and their inverses are Laurent polynomials over
a number field, and its standard invariant at one algebraic point
$\mathbf z_0$ is $\TL(\sqrt q)$. At every fixed depth the exceptional
phases lie in a proper algebraic subset of $(\C^\times)^d$. In all
degrees the conclusion therefore holds off a countable union of such
subsets, in particular at every algebraically independent tuple. For
$d=1$ the exceptional set in each degree is finite and algebraic.
\end{proposition}
\begin{proof}
The matrix $I-T_n(\mathbf z_0)$ has a nonzero minor of size
$q^n-\Cat_n$. The same minor is a Laurent polynomial in $\mathbf z$.
Where it is nonzero, the kernel has dimension at most $\Cat_n$.
It contains the canonical Temperley--Lieb space and hence equals it.
The assertions follow as in Theorem~\ref{thm:family}.
\end{proof}

For integers $r_1,\ldots,r_d$ the minor restricted to the curve
\[
 t\longmapsto H(z_{01}t^{r_1},\ldots,z_{0d}t^{r_d})
\]
is nonzero at $t=1$. The conclusion therefore holds at every
transcendental $t\in\T$ on this curve. The argument works as well
for algebraic-function entries on a finite cover, on each irreducible
component admitting a Temperley--Lieb specialization. On a normalized
curve, a limiting matrix suffices if its entries and inverses extend
there. The minor is then a nonzero rational function on that component.
Its zero set is finite and algebraic for data over $\overline\Q$.
In higher dimensions a finite morphism sends these exceptional
subvarieties into proper algebraic subsets of the parameter space.

For a monomial family with fourth-root coefficients in order
$q\equiv2\pmod4$, $q>2$, all tuples in $(\Q(i)\cap\T)^d$ give
Gaussian-rational entries. One can also transfer the local condition
from a matrix $G$ in Table~\ref{tab:ext-profiles}. Normalize
$v_p(p)=1$ and suppose $H$ has local unit entries with
\begin{equation}\label{eq:ext-close}
 v_p(H_{ay}-G_{ay})>v_p(\alpha)\qquad\text{for all }a,y,
\end{equation}
where we compare entries in a common local ring. The same strict
bound holds for the difference of the profiles, since
inversion preserves the valuation of a difference of units.
Consequently $S_H/\alpha$ and $S_G/\alpha$ have the same residue.
For a real base matrix the bound is $v_2(H_{ay}-G_{ay})>2$.
For a cube-root base matrix it is $v_3(H_{ay}-G_{ay})>1$.

In a monomial family recentered at the base point, it suffices to
impose this bound on each phase difference $z_j-1$. If the base
matrix is equivalent to the required representative, first apply that
fixed equivalence to the whole family. Between dephased Butson
representatives its row and column multipliers can be chosen among
roots of unity. The valuation bound then applies to the varying
phases $z_j$. For a real base, examples are
$z_j=(1+4it)/(1-4it)$, with $v_2(z_j-1)\geq3$.
For a cube-root base take $z_j=(1+9it\sqrt3)/(1-9it\sqrt3)$, with
$v_3(z_j-1)\geq5/2$. In both formulas $t\in\Z$ may be chosen
separately in each coordinate. The Gaussian-rational phases are
$(1+it)/(1-it)$, $t\in\Q$, together with $-1$.
\begin{table}[H]
\caption{Affine families to which Proposition~\ref{prop:ext-specialize}
applies. The last column gives additional algebraic parameters.
The family notation is that of \cite{LSO13,S08}.}\label{tab:ext-affine}
\centering\small
\begin{tabular}{@{}l>{\raggedright\arraybackslash}p{.42\textwidth}>{\raggedright\arraybackslash}p{.31\textwidth}@{}}\toprule
Order & Family & Algebraic phases covered\\ \midrule
$6$ & $D_6^{(1)}$ & Gaussian-rational phases and \eqref{eq:ext-haagerup-local} below\\[3pt]
$10$ & $D_{10}^{(3)},\ N_{10B}^{(3)},\ G_{10}^{(1)}$ & All Gaussian-rational tuples\\[3pt]
$14$ & $D_{14}^{(6)}$, including $D_{14}^{(5)}$ & All Gaussian-rational tuples\\[3pt]
$12$ & $H_{12B}^{(10)},\ H_{12C}^{(8)},\ B_{12D}^{(4)},\ B_{12E}^{(5)},$
$B_{12F}^{(7)},\ B_{12G}^{(7)},\ B_{12H}^{(6)},\ B_{12L}^{(5)}$
& Tuples satisfying \eqref{eq:ext-close} at a real specialization\\[3pt]
$4\bmod 8$, $q>4$ & Sz\"{o}ll\H{o}si's $(q/2+1)$-parameter families through a real matrix
& Tuples satisfying \eqref{eq:ext-close} at the real matrix\\ \bottomrule
\end{tabular}
\end{table}
In Table~\ref{tab:ext-affine} the order-ten and order-fourteen
formulas are fourth-root monomials \cite[Sections~5 and~7]{LSO13}.
All eight order-twelve families contain real class~1 of
\cite[Table~2]{LSO13}, whose equivalences also allow adjoint,
conjugation and transposition. The construction in the last row is
\cite[Theorem~3.5]{S08}. A Laurent family through a cube-root base
matrix as in Theorem~\ref{thm:ext-classes} has the same specialization
conclusion and the local bound at three.

\Needspace{13\baselineskip}
\subsection{Two families of order six}\label{ssec:ext-six}
The cube-root case in order six gives Tao's matrix \cite{T04}.
Haagerup's fourth-root matrix \cite[p.~307]{H97} belongs to
Di\c{t}\u{a}'s one-parameter family $D_6^{(1)}$
\cite{D04,TZ06}. We use the representative
\begin{equation}\label{eq:ext-affine}
 A(z)=\begin{pmatrix}
1&1&1&1&1&1\\
1&-1&i&i&-i&-i\\
1&i&-1&-i&z&-z\\
1&i&-i&-1&-z&z\\
1&-i&z^{-1}&-z^{-1}&i&-1\\
1&-i&-z^{-1}&z^{-1}&-1&i
\end{pmatrix},\qquad z\in\T.
\end{equation}

\begin{lemma}\label{lem:ext-haagerup}
The standard invariant of $A(z)$ is $\TL(\sqrt6)$ for every
transcendental phase. It is also $\TL(\sqrt6)$ for an algebraic
phase provided that at a prime above two one has
\begin{equation}\label{eq:ext-haagerup-local}
 v_2(z-\varepsilon)>\tfrac14\qquad
 \text{for some }\varepsilon\in\{1,-1,i,-i\},\qquad v_2(2)=1.
\end{equation}
This includes $z=(a+ib)/\sqrt{a^2+b^2}$ for integers $a,b$ of
opposite parity.
\end{lemma}
\begin{proof}
At $z=1$ the matrix is fourth-root, so
Proposition~\ref{prop:ext-specialize} gives the transcendental case.
For the local assertion consider the profile as a Laurent polynomial
over $\Z[i]$. Substituting the six rows of \eqref{eq:ext-affine} and
collecting powers of $z$ gives four possible residues
modulo $2\Z[i]$, namely
\[
 0,\quad (1+i)(1+z^2),\quad (1+i)(1+z^{-2}),\quad
 (1+i)(z+z^{-1}).
\]
In particular all coefficients are divisible by $1-i$, and the
derivative at each $\varepsilon$ in \eqref{eq:ext-haagerup-local}
is divisible by two. Write $s=v_2(z-\varepsilon)>1/4$. Since
$\varepsilon$ is a local unit, the terms of order at least two
in the expansion about $\varepsilon$ have valuation at least
$\tfrac12+2s$. The linear term has valuation at least $1+s$.
Thus $v_2(S_z-S_\varepsilon)>1$, and $S_z/2$ has the required
Gaussian residue.

For the last assertion take $a$ odd, $b$ even, and write
$c=\sqrt{a^2+b^2}$. The factors $c+a$ and $c-a$ differ by $2a$
of valuation one and have product $b^2$. If either factor had
valuation less than one, both would, contradicting $v_2(b^2)\geq2$.
Hence $v_2(z-1)\geq1$. For $b$ odd apply this to $-iz$.
\end{proof}

For Karlsson's two-parameter family \cite{K09}, the limits at the singular parameters
depend on the curve. We specify the phase coordinates and the square
roots used in taking those limits. Put
\[
 d=-\tfrac14(u-u^{-1})(v-v^{-1}),\qquad
 r_+^2=\frac{3-d}{1+d},\qquad r_-^2=\frac{3+d}{1-d},
\]
\[
 Z=\tfrac12\begin{pmatrix}1+u+v-uv&1-u+v+uv\\
                     1+u-v+uv&1-u-v-uv\end{pmatrix},\qquad
 \Sigma=\begin{pmatrix}r_+&r_-\\-r_-&-r_+\end{pmatrix}.
\]
Let $B(u)=\bigl(\begin{smallmatrix}1&1\\u&-u\end{smallmatrix}\bigr)$
and define $U_{jk}=-Z_{jk}(1+i\Sigma_{jk})/2$,
$V_{jk}=-Z_{jk}(1-i\Sigma_{jk})/2$. The matrix is
\begin{equation}\label{eq:ext-karlsson}
 K(u,v)=\begin{pmatrix}
 B(1)&B(u)&B(u)\\ B(v)^{\mathsf T}&U&V\\ B(v)^{\mathsf T}&V&U
 \end{pmatrix}.
\end{equation}
For $|d|<1$ we may choose the positive square roots, since changing
their signs gives row and column permutations \cite[Section~II]{K09}.

\begin{proposition}\label{prop:ext-karlsson}
The standard invariant of $K(u,v)$ is $\TL(\sqrt6)$ outside a
countable union of proper algebraic subsets of the parameter
space, in particular for algebraically independent $u,v\in\T$.
On the curves
\[
 (u,v)=(z,z^m)\quad(m\text{ positive and odd}),\qquad
 (u,v)=(z,z^{-1}),
\]
it is $\TL(\sqrt6)$ at every transcendental $z\in\T$ and at the
one-sided limits over $z=\pm i$. On the cubic curve the conclusion
also holds for algebraic $z\in\T$ with $v_2(z-i)\geq3$ at some
prime over two, including $z=(8-15i)/17$.
\end{proposition}

\begin{proof}
Consider first $z=e^{it}$, $v=z^m$ with $m$ positive and odd.
As $t\to\pi/2$, a shift by $\pi$ of the second angle if necessary
puts the angular displacements from the corner in the ratio $1:m$.
By \cite[Eqs.~(27)--(29)]{K09}, the limit is equivalent to
\eqref{eq:ext-affine}. Up to conjugation and a fourth-root factor
its parameter is
\[
 \frac{(1+m)+i(1-m)}{\sqrt{2(1+m^2)}}
 =\frac{a+ib}{\sqrt{a^2+b^2}},\qquad
 a=\frac{1+m}{2},\quad b=\frac{1-m}{2}.
\]
Since $a+b=1$, the
integers $a,b$ have opposite parity, and
Lemma~\ref{lem:ext-haagerup} applies. The sign and permutation
symmetries give the same conclusion for the remaining one-sided
limits over $\pm i$.

On the anti-diagonal $v=z^{-1}$ one has $d=-\sin^2t$,
$Z_{11}=-Z_{22}=\cos t$, $Z_{12}=1-i\sin t$, and
$Z_{21}=1+i\sin t$. As $t$ increases to $\pi/2$, $\cos t\,r_+\to2$
and $r_-\to1$. Substitution in \eqref{eq:ext-karlsson} gives a
fourth-root matrix, as do the other sides and sign choices.
Each component of the normalized cover, after completion, maps onto
the projective $z$-line. Its branches above $i$ therefore give
Temperley--Lieb specializations of the kind just computed.
Proposition~\ref{prop:ext-specialize} proves the assertion about
transcendental phases on each curve. The curve $m=1$ is the symmetric
Matolcsi--Sz\"{o}ll\H{o}si family. At $z=1$ the curves pass through
$F_6$ up to equivalence. There the condition for nonzero profile is
$a-c=b-d$ in $\Z/6\Z$, giving six components of the support of
$T_2$ and $\dim P_2=6$.

It remains to treat the algebraic phases on the cubic curve.
Write $\nu=v_2$. Take $z\ne i$ and suppose $s=\nu(z-i)\geq3$.
For $x=z+z^{-1}=(z-i)(z+i)/z$ we have $\nu(x)=s+1$.
Substituting $(u,v)=(z,z^3)$ in the formula for $d$ gives
\[
 d=-1+y,\qquad y=x^2(5-x^2)/4.
\]
The quantity $Y=xr_+$ is regular at the corner. It satisfies
\[
 Y^2=\frac{16}{5-x^2}-x^2,\qquad
 r_-^2=\frac{2+y}{2-y}.
\]
Choose boundary values $Y_0=4/\sqrt5$ and $r_{-,0}=1$. As
$\nu(y)=2s$, the preceding equations yield
\[
 \nu(Y^2/Y_0^2-1)\geq2s-2,\qquad \nu(r_-^2-1)=2s.
\]
For any unit $a$ with $\nu(a^2-1)>2$, the valuations of $a-1$
and $a+1$ are $1$ and $\nu(a^2-1)-1$ in some order. Apply this
to $Y/Y_0$ and $r_-$. After choosing signs,
$\nu(Y-Y_0)\geq2s-1$ and $\nu(r_--1)\geq2s-1$.

We compare the entries with their boundary values. Put
$J=1+z^2=zx$. The diagonal entries of $U,V$, multiplied by $-4$,
are of the form
\[
 (1\pm z-z^2)(J\pm izY).
\]
Here $\nu(J)=s+1$ and $\nu(zY-iY_0)\geq s+2$. Division by four
therefore gives differences of valuation at least $s-1$.
The corresponding off-diagonal products are
$(1\pm z\mp z^3+z^4)(1\pm ir_-)$. The polynomial changes by a term
of valuation at least $s$ and $\nu(1\pm i)=1/2$. After division
by four the valuation of the difference is at least
$\min\{s+1/2,2s-1\}-2>1$. In the first two rows and columns it
is at least $s$.

For $\eta=(2+i)/\sqrt5$, reordering the boundary rows by
$(1,2,4,6,5,3)$ and columns by $(1,2,3,5,4,6)$ gives
$A(-\bar\eta)$. Its normalized profile has the Gaussian residue
by the proof of Lemma~\ref{lem:ext-haagerup}. All corresponding
entry differences have valuation greater than one, so
\eqref{eq:ext-close} gives the same residue on the curve. This proves
the local assertion, including $z=i$ by the boundary calculation.
For the stated example use $(8-15i)/17-i=8(1-4i)/17$.

Finally the explicit pair $(z,z^3)$ at $z=(8-15i)/17$ has $|d|<1$.
Above this nonsingular pair all sign choices give equivalent
matrices. Over the open set $1-d^2\ne0$ the equations for $r_+$ and
$r_-$ define a finite cover. Each irreducible component dominates
this open set and therefore contains a point over the given pair.
Thus each component has a Temperley--Lieb specialization. Applying
Proposition~\ref{prop:ext-specialize} component by component and
adjoining the singular locus proves the assertion for the full
two-parameter family.
\end{proof}

\subsection{Butson-Hadamard examples}
Table~\ref{tab:ext-catalogue} lists the monomial equivalence classes
in the small Butson catalogue \cite{LOS20,ButsonData} to which
Theorem~\ref{thm:ext-classes} applies in the data considered here.
For each pair $(q,k)$ we used the exponent records in file order
(one matrix per line, read row by row) and numbered the
representatives starting at one. For each we checked $HH^*=qI_q$ and the indicated
profile congruence or sign reduction over $\Z$, $\Z[i]$ or
$\Z[\omega]$. The 752 fourth-root classes of order fourteen include
isolated matrices. Part~(2) of the theorem applies also to all
fourth-root matrices of order eighteen \cite{OP21} and, whenever
they exist, of every higher order congruent to two modulo four.
\begin{table}[htbp]
\caption{Catalogue representatives with Temperley--Lieb standard
invariant, grouped by order and allowed roots of unity.}\label{tab:ext-catalogue}
\centering\small
\begin{tabular}{@{}clr@{\hspace{1.8em}}l@{}}\toprule
Order & Class & Representatives & Profile calculation\\ \midrule
6 & $\mathrm{BH}(6,3)$ (Tao) & 1 & Third roots\\
6 & $\mathrm{BH}(6,4)$ (conference) & 1 & Gaussian rational\\
7 & $\mathrm{BH}(7,6)$ & 2 & Reflection modulo three\\
10 & $\mathrm{BH}(10,4)$ & 10 & Gaussian rational\\
12 & $\mathrm{BH}(12,2)$ & 1 & Real\\
12 & $\mathrm{BH}(12,3)$ & 2 & Third roots\\
14 & $\mathrm{BH}(14,4)$ & 752 & Gaussian rational\\
20 & $\mathrm{BH}(20,2)$ & 3 & Real\\
21 & $\mathrm{BH}(21,3)$ & 72 & Third roots\\ \midrule
 & Total & 844 & \\ \bottomrule
\end{tabular}
\end{table}

We give the sign reductions for the two order-seven matrices in the order
in which their representatives appear in $\mathrm{BH}(7,6)$ data file in the online catalogue of Butson
Hadamard matrices~\cite{ButsonData}. Reduce $w^e$ to $(-1)^e$.
In the first representative change the signs in row~3 and column~5.
In the second use row~5 and column~6. The unique negative entries
in successive rows then occur in columns $(5,6,1,3,2,7,4)$ and
$(6,5,2,3,1,7,4)$, respectively. Permuting the columns gives
$J_7-2I_7$ in both cases. The additional phases in the order-seven
family are described in Sections~\ref{sec:family}
and~\ref{sec:algebraic-family}.

For comparison, Jones's profile graph has $\dim P_2$ connected
components \cite[Theorem~2.11.15]{J21}. For each sixth-root
representative, form a vertex for every ordered pair $(a,c)$ and
join $(a,c)$ to $(b,d)$ when $S(a,b,c,d)\ne0$. This test uses the
two integer coefficients of $S$ in the basis $1,\omega$ of
$\Z[\omega]$. Counting components, including the diagonal block,
gives dimensions $6,4,4,2$ in $\mathrm{BH}(6,6)$ in file order.
The fourth representative belongs to Tao's class. In $\mathrm{BH}(8,6)$
nine of the 36 classes have dimension two, and in
$\mathrm{BH}(9,6)$ eight of the 17 do. The other classes in these
two orders have $\dim P_2>2$. All 436 classes in $\mathrm{BH}(13,6)$
have dimension two and all fail the reflection condition modulo
three. For these 436 classes and the indicated nine and eight,
the calculation leaves the higher relative commutants open.

\subsection{\textbf{Additional Examples}}\label{ssec:additional-examples}
In this section we give further examples in which our methods
determine the full Temperley--Lieb standard invariant. Table~\ref{tab:additional-examples}
gives the Jones index $I=[M:N]$ in each case. The invariant is
$\TL(\sqrt I)$ and both principal graphs are $A_\infty$.
For the Bisch--C\'aceres four-star our lowest-weight and specialization
arguments give a second proof. In the fixed-point examples the spin
results identify the acting quantum group.

\begingroup\small
\setlength{\LTpre}{5pt}\setlength{\LTpost}{0pt}
\makeatletter
\def\LT@makecaption#1#2#3{\LT@mcol\LT@cols{@{}c@{}}{%
  \parbox{\linewidth}{\@tempcnta=64\relax\@makecaption{#1{#2}}{#3}}}}
\makeatother
\setlength{\baselineskip}{10.4pt}
\setlength{\tabcolsep}{3.5pt}
\renewcommand{\arraystretch}{1.02}
\begin{longtable}{@{}>{\raggedright\arraybackslash}p{.232\linewidth}%
                    >{\raggedright\arraybackslash}p{.101\linewidth}%
                    >{\raggedright\arraybackslash}p{.278\linewidth}%
                    >{\raggedright\arraybackslash}p{\dimexpr.389\linewidth-6\tabcolsep\relax}@{}}
\caption{Hyperfinite subfactors with Temperley--Lieb standard invariant.}\label{tab:additional-examples}\\
\toprule
\textbf{Class of examples} & \textbf{Index $I$} & \textbf{Proof of TL} &
\textbf{Relation to earlier results}\\
\midrule
\endfirsthead
\multicolumn{4}{c}{\textsc{Table \thetable} (continued)}\\[3pt]
\toprule
\textbf{Class of examples} & \textbf{Index $I$} & \textbf{Proof of TL} &
\textbf{Relation to earlier results}\\
\midrule
\endhead
\multicolumn{4}{@{}l}{\textbf{Spin models (Hadamard commuting squares)}}\\*[2pt]
Two-entry Hadamards with nonreal ratio $z$ (Fano plane,
Bj\"orck primes, Sylvester cores).
& $4m-1$ for each such design, e.g.\ $7,11,15$
& At a prime $p\mid I$ the design identity $(z-1)^2/z=-I/m$ gives
$2v(z-1)=v(I)$. This yields the four-row congruence.
& Outside Theorem~\ref{thm:ext-classes}. The Fano matrix lies outside
Petrescu's family up to equivalence. For the designs see \cite{Sz11}.
We know of no earlier full TL proof for this class.\\[3pt]
Beauchamp--Nicoara self-adjoint family, for the admissible $y$ at right.
& $6$
& For transcendental $y$, specialize at $y=i$. For algebraic $y$ with
$v_2(y-i)>3$ at some place, transfer the profile of the fourth-root
matrix at $i$.
& Proposition~\ref{prop:ext-specialize} and \eqref{eq:ext-close}
give these cases of the conjecture in \cite{BeN08}, including $i$.
Here $v_2(2)=1$. We know of no independent proof for these ranges.\\
\midrule
\multicolumn{4}{@{}l}{\textbf{Non-spin commuting squares}}\\*[2pt]
Stojanovic's four infinite-dimensional multi-matrix squares.
& $16/3$, $25/4$, $36/5$, $64/7$
& Counting $2^n$ paths far down the tail gives $\dim P_n\le4^n$
for odd $n$, hence graph norm at most two as in Section~\ref{sec:annular}.
No $A_\infty$ tail result applies.
& Stojanovic \cite{St21} gave the inclusions and indices.
We apply \cite[Theorem~17.4]{Sch90}.
We know of no earlier full TL proof for these four inclusions.\\[3pt]
Bisch--C\'aceres four-star $S(3,3,3,3)$ at $w=-1$ and all but
finitely many admissible phases on the same branch.
& $3+\sqrt5$
& At $w=-1$ leaf path counts give TL through depth $3$.
Flatness excludes extra lowest-weight vectors at depths $4,5$.
The classification \cite{AMP} then gives TL at all depths.
& Bisch--C\'aceres proved existence by countability and an intermediate
subfactor obstruction \cite[Theorem~6.6]{BC25}. Here the flatness
minors at $w=-1$ also give the assertion on the branch by specialization.\\[3pt]
Four-star connection on $S(1,3,3,3)$.
& $5$
& Flatness rank $2$ at depth $4$ modulo $3121$.
& Schou proved existence of a biunitary on this graph \cite[\S3]{Sch90}.\\[3pt]
Five-star connection on $S(1,1,1,1,2)$.
& $3+\sqrt5$
& Flatness ranks $3,8,24$ at depths $3,4,5$ modulo $3121$.
& Outside the four-star families of \cite{BC25}.\\
\midrule
\multicolumn{4}{@{}l}{\textbf{Quantum-group fixed points (spin-square realization for the canonical action)}}\\*[2pt]
Mixed model $(H,K)$ with one TL Hadamard input, e.g.\ Tao and $F_6$.
& Common order $q>4$
& Compression recovers the single-matrix model. Its Hopf image
forces the mixed Hopf image to be $S_q^+$.
& Banica computed the $S_q^+$ fixed-point invariant \cite{Ba99}.
The TL input identifies the group for the mixed model of \cite{BaN17}.\\[3pt]
Universal flat model at orders six and seven.
& $6,\ 7$
& Evaluation at Tao's matrix or Petrescu's matrix gives the full
Hopf image $S_q^+$.
& This follows from the spin results above.
The independent inner-faithfulness theorem in
\cite[Corollary~4.9]{BCF20} covers $q\le5$ and $q\ge10$.\\
\bottomrule
\end{longtable}
\par\vspace{5pt}
\noindent\begin{minipage}{\linewidth}
The last two rows give $P^G\subset(\C^q\otimes P)^G$ with $G=S_q^+$.
For the canonical action attached to $H$, this is the spin-square subfactor
of $H$ \cite{Ba00}.
\end{minipage}
\endgroup

\appendix

\section{Equivalence of Petrescu's matrix $P_7$ and the matrix $H$
of Section~\ref{sec:matrix}}\label{app:equivalence}
In \cite[Eq.~(96)]{TZ06} the representative $P_7=(w^{E_{ij}})$
has exponent matrix
\[
E=\begin{pmatrix}
0&0&0&0&0&0&0\\
0&1&4&5&3&3&1\\
0&4&1&3&5&3&1\\
0&5&3&1&4&1&3\\
0&3&5&4&1&1&3\\
0&3&3&1&1&4&5\\
0&1&1&3&3&5&4
\end{pmatrix}.
\]
Applying the row and column operations of Lemma~\ref{lem:matrix}
gives the exponent matrix
\[
\begin{pmatrix}
3&0&0&0&0&0&0\\
0&1&0&4&2&4&0\\
0&0&1&0&4&2&4\\
0&4&0&1&0&4&2\\
0&2&4&0&1&0&4\\
0&4&2&4&0&1&0\\
0&0&4&2&4&0&1
\end{pmatrix},
\]
whose first row and circulant core give \eqref{eq:H}. All entries
and operations lie in $\Z[w]$, with $w^2=w-1$.

\end{document}